\documentstyle{amltd}
\begin{document}
\annalsline{151}{2000}
\received{July 8, 1997}
 \startingpage{59}
\def\bsp{\begin{eqnarray} \noalign{\vskip4pt}
&&}
\def\overset#1#2{\stackrel{#1}{#2}}
 
\catcode`\@=11
\font\twelvemsb=msbm10 scaled 1100
\font\tenmsb=msbm10
\font\ninemsb=msbm10 scaled 800
\newfam\msbfam
\textfont\msbfam=\twelvemsb  \scriptfont\msbfam=\ninemsb
  \scriptscriptfont\msbfam=\ninemsb
\def\msb@{\hexnumber@\msbfam}
\def\Bbb{\relax\ifmmode\let\next\Bbb@\else
 \def\next{\errmessage{Use \string\Bbb\space only in math
mode}}\fi\next}
\def\Bbb@#1{{\Bbb@@{#1}}}
\def\Bbb@@#1{\fam\msbfam#1}
\catcode`\@=12

 \catcode`\@=11
\font\twelveeuf=eufm10 scaled 1100
\font\teneuf=eufm10
\font\nineeuf=eufm7 scaled 1100
\newfam\euffam
\textfont\euffam=\twelveeuf  \scriptfont\euffam=\teneuf
  \scriptscriptfont\euffam=\nineeuf
\def\euf@{\hexnumber@\euffam}
\def\frak{\relax\ifmmode\let\next\frak@\else
 \def\next{\errmessage{Use \string\frak\space only in math
mode}}\fi\next}
\def\frak@#1{{\frak@@{#1}}}
\def\frak@@#1{\fam\euffam#1}
\catcode`\@=12

\font\emi= cmmi10 scaled 1700 
\font\eightmi=cmmi10

\font\tenrm=cmr10
\def\bye{\end{document}}
\title{Integral mappings and the principle of\\ local
reflexivity for
noncommutative\\ {\emi L}\raise8pt\hbox{\eightmi 1}-spaces}
\shorttitle{Local reflexivity} 
\acknowledgements{The research of E. Effros and Z.-J. Ruan was partially supported by the National Science
Foundation.\hfill\break
\hglue.25in 1991 {\it Mathematics Subject Classification\/}: Primary 47D15 and 46B07; Secondary 46B08.}
 \twoauthors{Edward G. Effros, Marius Junge,}{Zhong-Jin Ruan}
   \institutions{UCLA, Los Angeles, CA\\
{\eightpoint {\it E-mail address\/}: ege@math.ucla.edu}\\
\vglue6pt
University of Illinois, Urbana, IL\\
{\eightpoint {\it E-mail address\/}: junge@math.uiuc.edu} \\ \vglue6pt
 University of Illinois, Urbana, IL\\
{\eightpoint {\it E-mail address\/}: ruan@math.uiuc.edu}}

\centerline{\it Dedicated to the memory of Irving Segal}

\bigbreak

\centerline{\bf Abstract}
\bigbreak

The operator space analogue of the {\em strong form} of the principle of
local reflexivity
is shown to hold for any von Neumann algebra predual, and thus for any
$C^{*}$-algebraic dual. This is in striking contrast to the situation for
$C^{*}$-algebras, since, for example, $K(H)$ does not have that property. The
proof uses the Kaplansky density theorem together with a careful analysis
of two
notions of integrality for mappings of operator spaces.

\section{Introduction}

Transcendental models, such as ultraproducts and second duals of
nonreflexive spaces, arise quite naturally in Banach space theory. Despite
their esoteric nature, these constructions have proved to be indispensable
for the classification of von Neumann algebras and $C^{*}$-algebras (see,
e.g.,  \cite{Con}, \cite{Lo}, and \cite{McD}). Generally speaking, if one
wishes to prove that a Banach space or a $C^{*}$-algebra has some
approximate property, one begins by proving that an appropriate model has
the corresponding exact property. One must then relate the exact property in
the model to the approximate property in the original space. In Banach space
theory, this is often accomplished by using the {\it principle of local
reflexivity}.

In its weakest form, which was first proved by Schatten in his early
monograph \cite{S}, the principle of local reflexivity states that any
finite-dimensional subspace $L$ of the second dual $E^{**}$ of a Banach
space $E$ can be approximated by finite-dimensional subspaces of $E$ in the
weak$^{*}$ topology. The importance of this result became evident in
Grothendieck's ground-breaking study of Banach spaces in the 1950's (see,
e.g., \cite{G}, \cite{Pi}). His theory rested, in part, upon relating
various canonical tensor products to corresponding mapping spaces. One of
his key observations, which is equivalent to the principle of local
reflexivity, is that if $E$ and $F$ are Banach spaces, then
\begin{equation}
(E\otimes _{\lambda }F)^{*}=I(E,F^{*}),  \label{classical}
\end{equation}  
where $\otimes _{\lambda }$ is the injective Banach space tensor product, and $%
I(E,F^{*})$ is the space of integral mappings $\varphi :E\rightarrow F^{*}.$
In the 1960's, Johnson, Lindenstrauss, Rosenthal, and Zippin (see \cite{D}, \cite{Jo},
 and  \cite{Li}) formulated a strong form of local reflexivity,
which implies that the approximating subspaces of $E$ are close to $L$ in
the Banach-Mazur distance.

In recent years it has become evident that one can adapt Banach space
techniques to the study of linear spaces of Hilbert space operators,
provided one replaces the  bounded linear mappings of Banach space theory by the
{\it completely bounded} linear mappings (see \cite{R}). As a result,
there has been a remarkable convergence of classical and ``noncommutative''
functional analysis. Much of operator space theory has been developed along
the lines pioneered by Grothendieck (see, e.g., \cite{BP}, \cite{ER1}, \cite
{ER3} and \cite{ER4}).

The operator space analogue of the weak form of local reflexivity was
introduced in \cite{EH}, and was further studied in \cite{ER1}, \cite{ER3},
\cite{ER4}, \cite{ER5} and \cite{J}. An operator space $V$ is defined to be
{\it locally reflexive} if for each finite-dimensional operator space $F,$
any complete contraction $\varphi :F\rightarrow V^{**}$ may be approximated
in the point-weak$^{*}$ topology by a net of complete contractions $\varphi
_{\lambda }:F\rightarrow V.$ Equivalently, $V$ is locally reflexive if and
only if for every finite-dimensional operator space $F,$ we have the
isometry
\begin{equation}
F^{*}\check{\otimes}V^{**}=(F^{*}\check{\otimes}V)^{**}  \label{F1.0}
\end{equation}
(this is essentially condition $C^{\prime \prime }$ introduced in \cite{AB},
\cite{EH}), or what is the same, $V$ is locally reflexive if and only if we
have the natural isometry
\begin{equation}
(F^{*}\check{\otimes}V)^{*}=F\hat{\otimes}V^{*}  \label{F1.1}
\end{equation}
for each finite-dimensional operator space $F$. All exact $C^{*}$-algebras
are locally reflexive (see \cite{Ki1}). On the other hand, it was shown in
\cite{EH} that some $C^{*}$-algebras are {\it not} locally reflexive. The
strong version of local reflexivity does not seem to have an interesting $%
C^{*}$-algebraic analogue, since apparently few $C^{*}$-algebras have that
property (see \S 6).

Turning to other operator spaces, the second author showed that the operator
space $T(H)$ of trace class operators on a Hilbert space $H$ is locally
reflexive \cite{J}. The argument is unexpectedly subtle. The proof used
asymptotic techniques related to Pisier's ultraproduct theory (see \cite{Pi1}%
), as well as a novel application of the Kaplansky density theorem (see \S 2
and \S 7). Employing different methods, the first and third author extended
this result to the preduals of injective von Neumann algebras \cite{ER5}.

In this paper, we prove that the predual of {\it any} von Neumann algebra
is locally reflexive. We recall that the space $T(H)$ may be regarded as a
``noncommutative $\ell ^{1}$-space'', and in turn, the preduals of von
Neumann algebras play the role of the ``noncommutative $L^{1}$-spaces''
mentioned in the title of this paper. What is even more surprising is that
these operator spaces are locally reflexive in the {\it strong sense, }%
i.e., we can assume that the approximations are close in the sense of the
Pisier-Banach-Mazur distance for operator spaces. The approach in this paper
is rather different than that used in either \cite{ER5} or \cite{J}, since
it does not depend upon ultraproduct techniques.

As in \cite{ER5} and \cite{J}, the Kaplansky density theorem plays a
fundamental role in this paper. We begin in Section~2 by showing how that result
implies an unexpected theorem about completely bounded mappings $\varphi
:A^{*}\rightarrow B$ for $C^{*}$-algebras $A$ and $B.$ Our analysis of local
reflexivity rests upon a careful study of the relationship between the
completely nuclear, completely integral, and exactly integral mappings
introduced in \cite{ER3}, \cite{ER4}, and \cite{J}, respectively. These
results are presented in Section~3, Section~4, and Section~5, respectively. The notion of
exactly integral mappings is the most novel of these definitions, and we have
explored it in considerable detail in Section~5. As we have indicated in the text, much
of the material in Section~5 is not needed in the subsequent sections.

The main theorem on local reflexivity is proved in Section~6 (Theorem \ref{P5.7}). In Section~7 we apply Theorem
\ref{P5.7} to show that the preduals of von Neumann algebras with the {\rm QWEP} property of Kirchberg and Lance (see
\cite {Ki1}, \cite{La}) are locally approximable by subspaces of dual matrix
spaces $T_{n}\ $with $n\in \Bbb{N}$ (see below). This covers a remarkably
large class of von Neumann algebras, and in fact it has been conjectured
that {\it all }$C^{*}$-algebras have the {\rm QWEP}. We conclude by showing
that the main theorem also implies a factorization theorem that was used by
the second author in his proof that $T(H)$ is locally reflexive (see above).

Given any Hilbert space $H,$ we let $B(H),$ $T(H),$ and $K(H)$ denote the
bounded, trace class, and compact operators on a Hilbert space $H,$ and we
let $M_{n},$ $T_{n},$ and $K_{n}$ denote these operator spaces when $H=\Bbb{%
\Bbb{C}}^{n}$ for $n<\infty $ and $H=l^{2}$ for $n=\infty .$ We use the
pairings
\begin{equation}
\langle a,b\rangle =\sum a_{i,j}b_{i,j}  \label{Duality}
\end{equation}
for $a\in K_{\infty }$ or $M_{\infty },$ and $b\in T_{\infty }.$ Given
operator spaces $V$ and $W,$ we\break let ${\cal  C}{\cal  B}(V,W)$ denote the
operator space of completely bounded mappings\break $\varphi :V\rightarrow W,$
with the norm
$$
\left\| \varphi \right\| _{cb}=\sup \left\{ \left\| {\rm id}\otimes \varphi
:M_{n}\otimes V\rightarrow M_{n}\otimes W\right\| \right\} .
$$
If $V$ and $W$ are operator spaces, we have corresponding injective and
projective operator space tensor products $V\check{\otimes}W$ and $V\hat{%
\otimes}W$. For the first, let us suppose that we have concrete
representations $V\subseteq {\cal  B}(H)$ and $W\subseteq {\cal  B}(K).$
Then $V\check{\otimes}W$ is defined to be the closure of $V\otimes W$ in $%
{\cal  B}(H\otimes K).$ On the other hand, the operator space $V\hat{%
\otimes}W$ is uniquely determined by the fact that we have a complete
isometry:
$$
{\cal  (}V\hat{\otimes}W)^{*}\cong {\cal  CB}(V,W^{*}).
$$
We write $V\otimes _{\vee }W$ and $V\otimes _{\wedge }W$ for the algebraic
tensor product together with the relative matrix norms.

We emphasize that although we have used Banach space notation for these
tensor products, they generally do not coincide with the corresponding
tensor products of Banach space theory$.$ On the other hand, the properties
of these operator space tensor  products under mappings are quite analogous to
their Banach space antecedents (see, e.g., \cite{BP} and \cite{ER1}). We
also appropriate the Banach space notation $\nu $ and $\iota $ for the {\it %
\ completely} nuclear and {\it completely} integral mapping norms (see \S 3
and \S 4).

We shall say that an operator space is a {\it matrix space} if it is
completely isometric to a subspace of $M_{n}$ for some $n\in \Bbb{N}.$
Unless otherwise indicated, we consider only complete operator spaces. For
our purposes it often suffices to regard various mapping spaces as Banach
spaces rather than operator spaces, i.e., we do not need to consider the
natural matrix norms on these spaces. Reflecting this, we will at times
state that a mapping is a ``(complete) contraction,'' or a ``(complete)
quotient mapping'' to indicate that although it is true, there is no need to
prove the stronger statement.

In order to make this paper more accessible to operator algebraists, we have
largely avoided using the formal machinery of operator ideals. Given a pair
of operator spaces, we identify the algebraic tensor product $V^{*}\otimes W$
with the vector space ${\cal  F}(V,W)$ of continuous finite rank mappings $%
\varphi :V\rightarrow W.$ We use the terminology ``operator ideal'' to mean
an assignment to each pair of operator spaces $V$ and $W,$ a space of
completely bounded mappings $\alpha (V,W)\supseteq {\cal  F}(V,W),$ with a
norm $\alpha (\varphi )$, such that
\begin{equation}
\alpha (\tau \circ \varphi \circ \sigma )\le \left\| \tau \right\|
_{cb}\alpha (\varphi )\left\| \sigma \right\| _{cb}  \label{opid}
\end{equation}
whenever we are given a diagram of mappings
$$
X\overset{\sigma }{\longrightarrow }V\overset{\varphi }{\longrightarrow }W%
\overset{\tau }{\longrightarrow }Y.
$$
We let $\alpha ^{0}(V,W)$ denote ${\cal  F}(V,W)$ with the relative norm
in $\alpha (V,W)$.

Given a Banach space $V,$ we have a corresponding linear mapping
$$
\hbox{{\rm trace}}:{\cal  F}(V,V)=V^{*}\otimes V\rightarrow \Bbb{C}
$$
defined by
$$
\hbox{{\rm trace\thinspace }}(f\otimes x)=f(x)
$$
for $x\in V$ and $f\in V^{*}$. Given Banach spaces $V$ and $W$ and bounded
linear mappings $\varphi :V\rightarrow W$ and $\psi :W\rightarrow V,$ with $%
\psi \in {\cal  F}(W,V),$ we have the corresponding {\it trace duality
pairing}
\begin{equation}
\langle \varphi ,\psi \rangle =\hbox{{\rm trace\thinspace }}(\varphi
\circ \psi )=\hbox{{\rm trace\thinspace }}(\psi \circ \varphi ).
\label{traceduality}
\end{equation}
If we let
$$
\psi =\sum_{i=1}^{n}g_{i}\otimes v_{i}\in W^{*}\otimes V,
$$
we have
\begin{equation}
\langle \varphi ,\psi \rangle =\hbox{{\rm trace}}\left(\sum_{i}g_{i}\otimes
\varphi (v_{i})\right)=\hbox{{\rm trace}}\,({\rm id}\otimes \varphi )(\psi ).
\label{traceduality2}
\end{equation}
Finally we note that for any operator space $V,$
\begin{equation}
\hbox{{\rm trace: }}V^{*}\otimes _{\wedge }V\rightarrow \Bbb{C}
\label{traceduality3}
\end{equation}
is contractive since $(f,v)\mapsto f(v)$ is a completely contractive bilinear mapping.

\section{Finite rank approximations and the\\ Kaplansky density theorem}

Given operator spaces $V$ and $W,$ we say that a completely bounded mapping $%
\varphi :V^{*}\rightarrow W$ satisfies the {\it weak}$^{*}$ {\it %
approximation property} (${\rm W}^{*}{\rm AP}$) if there exists a net of finite rank weak%
$^{*}$ continuous mappings $\varphi _{\lambda }:V^{*}\rightarrow W$ with $%
\left\| \varphi _{\lambda }\right\| _{cb}\le \left\| \varphi \right\| _{cb}$
which converge to $\varphi $ in the {\it point-norm} topology. If $H$ is an
infinite-dimensional Hilbert space, the identity mapping $I:B(H)\rightarrow
B(H)$ does not have such approximations since $B(H)$ does not have the
metric approximation property of Grothendieck \cite{Sz}. Our object in this
section is to show that by contrast, if $A$ and $B$ are $C^{*}$-algebras,
{\it any }completely bounded mapping $\varphi :A^{*}\rightarrow B$ has the $%
{\rm W}^{*}{\rm AP}.$

Given von Neumann algebras $R$ and $S$, we let $R\overline{\otimes }S$
denote the von Neumann algebra tensor product of $R$ and $S$. Then each
function $f\in R_{*}$ determines a {\it slice mapping}
$$
f\otimes {\rm id}:R\overline{\otimes }S\rightarrow S
$$
where
$$
\left\langle (f\otimes {\rm id})(u),~g\right\rangle =(f\otimes g)(u)
$$
for $u\in R\overline{\otimes }S$ and $g\in S_{*}$ (see \cite{T}). As a
result, each element $u$ in $R\overline{\otimes }S$ determines a mapping $%
\varphi _{u}\in {\cal  C}{\cal  B}(R_{*},S)$ by
$$
\varphi _{u}(f)=(f\otimes {\rm id})(u).
$$
It was shown in \cite{ER1} that this determines a complete isometry
\begin{equation}
R\overline{\otimes}S\cong {\cal  CB}(R_{*},S).  \label{2.1}
\end{equation}

\proclaim{Lemma}\label{P2.1}
Given $C^{*}$\/{\rm -}\/algebras $A$ and $B${\rm ,} every complete contraction $
\varphi :A^{*}\rightarrow B^{**}$ can be approximated by a net of finite
rank weak$^{*}$ continuous complete contractions $\varphi _{\lambda
}:A^{*}\rightarrow B$ in the point\/{\rm -}\/weak$^{*}$ topology{\rm .} Every completely
bounded mapping $\varphi :A^{*}\rightarrow B$ satisfies the ${\rm W}^{*}{\rm AP}.$
\endproclaim

\demo{Proof}
Using the universal representations of $A$ and $B$, we may identify $A^{**}$
and $B^{**}$ with von Neumann algebras on Hilbert spaces $H$ and $K$. The\break $%
^{*}$-algebra $A\otimes B$ is weak operator dense in the von Neumann algebra
$A^{**}\overline{\otimes }B^{**}=(A\otimes B)^{\prime \prime }$ on $H\otimes
K$. From the Kaplansky density theorem, the unit ball of the $^{*}$-algebra $%
A\otimes B$ is weak operator dense in that of $A^{**}\overline{\otimes }%
B^{**}$.

If $\varphi :A^{*}\rightarrow B^{**}$ is a complete contraction, we may
assume that $\varphi =\varphi _{u}$ for some contractive element $u\in A^{**}%
\overline{\otimes}B^{**}$ since we have the isometry
$$
A^{**}\overline{\otimes}B^{**}={\cal  CB}(A^{*},B^{**})
$$
by (\ref{2.1}). There exists a net of contractive elements $u_{\lambda }\in
A\otimes B$ converging to $u$ in the weak operator topology on $B(H\otimes
K) $. It follows that $u_{\lambda }$ converges to $u$ relative to the
topology determined by the algebraic tensor product $A^{*}{\otimes }B^{*}$.
We have that $\varphi _{\lambda }=\varphi _{u_{\lambda }}$ is a net of
finite rank weak$^{*}$ continuous complete contractions from $A^{*}$ into $B$
which converges to $\varphi =\varphi _{u}$ in the point-weak$^{*}$ topology,
i.e., for each $f\in A^{*},$
$$
\varphi _{\lambda }(f)=f\otimes id(u_{\lambda })\in B\rightarrow \varphi
(f)=f\otimes id(u)\in B^{**}
$$
in the weak$^{*}$ topology.

If $\varphi $ is a complete contraction from $A^{*}$ into $B$, we have that $%
\varphi _{\lambda }$ converges to $\varphi $ in the point-weak topology.
The usual convexity argument shows that we can find a net of finite rank weak%
$^{*}$ continuous complete contractions\break $\psi _{\mu}:A^{*}\to B$
in the convex hull of $\{\varphi _{\lambda }\}$, which converges to $\varphi
$ in the point-norm topology (see, e.g., \cite[p.~477]{DS}).
\enddemo  
\vglue-12truept

\section{Completely nuclear mappings}
\vglue-6truept
Given operator spaces $V$ and $W,$ there is a canonical mapping
\begin{equation}
\Phi :V^{*}\hat{\otimes}W\rightarrow V^{*}\check{\otimes}W\subseteq {\cal  %
CB}(V,W)  \label{F4}
\end{equation}
which extends the identity mapping on the algebraic tensor product $%
V^{*}\otimes W$. A linear mapping $\varphi :V\rightarrow W$ is said to be
{\it completely nuclear} if it lies in the image of $\Phi $ (see \cite{ER3}). Identifying the linear space ${\cal 
N}(V,W)=\Phi (V^{*}\hat{\otimes}W)$ with the quotient Banach space $V^{*}\hat{\otimes}W/{\rm ker\,\,}\Phi ,$
we call the corresponding norm $\nu $ the {\it completely nuclear norm}  
 on ${\cal  N}(V,W).$ If $V$ or $W$ is
finite-dimensional, then $\Phi
$ is one-to-one, and we have the isometry
\begin{equation}
{\cal  N}(V,W)=V^{*}\hat{\otimes}W.  \label{nuclearfin}
\end{equation}
If $\varphi :V\rightarrow W$ is a linear mapping which is not nuclear, we
write $\nu (\varphi )=\infty .$

Turning to a ``prototypical'' example, let us suppose that $a$ and $b$ are
infinite scalar matrices with Hilbert-Schmidt norms $\left\| a\right\|
_{2},\left\| b\right\| _{2}<1.$ Then the mapping
$$
M(a,b):M_{\infty }\rightarrow T_{\infty }:x\rightarrow axb,
$$
satisfies $\nu (M(a,b))<1.$ More generally, given any operator spaces $V$
and $W,$ a linear mapping $\varphi :V\rightarrow W$ satisfies $\nu (\varphi
)<1$ if and only if it factors through such a mapping via completely
contractive mappings. Thus we have that $\nu (\varphi )<1$ if and only if
there is a commutative diagram
\begin{eqnarray}
&& \label{F2}\\
\noalign{\vskip-16pt}
&&\begin{array}{ccc}
M_{\infty } & \overset{M(a,b)}{\longrightarrow } & T_{\infty } \\
{\scriptstyle r}\uparrow &  & \downarrow {\scriptstyle s} \\
V & \overset{\varphi }{\longrightarrow } & W
\end{array}
\nonumber
\end{eqnarray}  
where $r$ and $s$ are complete contractions, and $\left\| a\right\|
_{2},\left\| b\right\| _{2}<1$. It is also equivalent to assume that there
is a commuting diagram
\begin{eqnarray}
&& \label{F3}\\
\noalign{\vskip-16pt}
&&
\begin{array}{ccc}
K_{\infty } & \overset{M(a,b)}{\longrightarrow } & T_{\infty } \\
{\scriptstyle r}\uparrow &  & \downarrow {\scriptstyle s} \\
V & \overset{\varphi}{\longrightarrow } & W
\end{array}
\nonumber
\end{eqnarray}
with the same assumptions (see \cite{ER3}).  

\proclaim{Lemma}
\label{P3.0} 
Given operator spaces $V$ and $W${\rm ,} the canonical mapping
$$
{\rm id}\otimes \iota _{W}:V\hat{\otimes}W\rightarrow V\hat{\otimes}W^{**}
$$
is a complete isometry{\rm .}
\endproclaim

 \demo{Proof} 
Let $\iota _{W}:W\rightarrow W^{**}$ be the canonical embedding. It follows
from the definition of the projective tensor product that ${\rm id}\otimes \iota _{W}$ is
a complete contraction from $V\hat{\otimes}W$ into $V\hat{\otimes}W^{**}$.
In order to show that ${\rm id}\otimes \iota _{W}$ is isometric, it suffices to
show that its adjoint $({\rm id}\otimes \iota _{W})^{*}$ is a norm quotient
mapping. Equivalently, since we have the commutative diagram
$$
\begin{array}{ccc}
(V\hat{\otimes}W^{**})^{*} & \overset{({\rm id}\otimes \iota _{W})^{*}}{%
\longrightarrow } & (V\hat{\otimes}W)^{*} \\
|| &  & || \\
{\cal  {\cal  C}B}(V,W^{***}) & \overset{\theta }{\longrightarrow } &
{\cal  CB}(V,W^{*}) \ ,
\end{array}
$$
where $\theta (\varphi )=(\iota _{W})^{*}\circ \varphi ,$ it suffices to
prove that $\theta $ is a quotient mapping. If we are given a complete
contraction $\psi :V\rightarrow W^{*},$ we have that
$$
\psi =(\iota _{W})^{*}(\iota _{W^{*}}\circ \psi ),
$$
where $\iota _{W^{*}}\circ \psi $ is a complete contraction in ${\cal  C}%
{\cal  B}(V,W^{***}).$ Thus $\theta $ is indeed a quotient mapping, and $%
{\rm id}\otimes \iota _{W}$ is isometric.

Applying this to the space $T_{n}\hat{\otimes}V,$ and using associativity of
the projective tensor product, it follows that we have an isometry
\begin{eqnarray*}
&&\\
\noalign{\vskip-9pt}
&&
{\rm id}\otimes ({\rm id}\otimes \iota _{W}):T_{n}\hat{\otimes}(V\hat{\otimes}W)\to T_{n}%
\hat{\otimes}(V\hat{\otimes}W^{**}) \\
\noalign{\vskip-9pt}
\end{eqnarray*}
for each $n\in \Bbb{N}$. Taking the adjoint,
\begin{eqnarray*}
&&\\
\noalign{\vskip-9pt}
&&
({\rm id}\otimes \iota _{W})_{n}^{*}:M_{n}((V\hat{\otimes}W^{**})^{*})\rightarrow
M_{n}((V\hat{\otimes}W)^{*})
\\
\noalign{\vskip-9pt}
\end{eqnarray*}
is a quotient mapping, and thus $({\rm id}\otimes \iota _{W})^{*}$ is a complete
quotient mapping. It follows that $({\rm id}\otimes \iota _{W})^{**}$ is a
complete isometry, and restricting it to $V\hat{\otimes}W,$ we conclude that
${\rm id}\otimes \iota _{W}$ is a complete isometry.%
\enddemo %

\vglue6pt
\proclaim{Lemma}
\label{P3.1} Given a nuclear mapping $\varphi :V\rightarrow W,$ we have $\nu
(\varphi ^{*})\le \nu (\varphi ).$ If $V$ or $W$ is finite-dimensional{\rm ,} then
$\nu (\varphi ^{*})=\nu (\varphi ).$
\endproclaim
\vglue6pt

\demo{Proof}%
If we let $S(\varphi )=\varphi ^{*},$ it is evident from the commutative
diagram
\begin{eqnarray}
&&\nonumber\\
\noalign{\vskip-9pt} 
&& \label{diag}
\\
\noalign{\vskip-16pt}
&&\begin{array}{ccc}
V^{*}\hat{\otimes}W & \overset{{\rm id}\otimes \iota _{W}}{\longrightarrow } &
V^{*}\hat{\otimes}W^{**} \\
\downarrow {\scriptstyle \Phi _{1}} &  & \downarrow {\scriptstyle \Phi _{2}}
\\
{\cal  N}(V,W) & \overset{S}{\longrightarrow } & {\cal  N}(W^{*},V^{*})
\end{array}
\nonumber\\
\nonumber
\end{eqnarray}
that $S$ is a contraction. Even though the top row is isometric (Lemma \ref
{P3.0}), and the two columns are quotient mappings, it does not follow that $%
S$ is isometric, since one might have that
\begin{eqnarray*}
&&\\
\noalign{\vskip-9pt}
&&
\ker \Phi _{2}\cap (V^{*}\hat{\otimes}W)\neq \ker \Phi _{1}.
\\
\noalign{\vskip-9pt}
\end{eqnarray*}
On the other hand, if either $V$ or $W$ is finite-dimensional, then
the mappings $\Phi
_{i}$ are isometric, and thus the same is true for $S.$%
\enddemo %

\vglue6pt

We note that if $V$ and $W$ are both infinite-dimensional, we can have that $
\nu (\varphi ^{*})<\nu (\varphi )$ even if $\varphi $ is of finite rank (see
\cite[p.~67]{DF}). On the other hand if $V^{*}$ has the operator
approximation property (see \cite{ER1}), the mappings $\Phi _{i}$ in (\ref
{diag}) are one-to-one, and thus $\nu (\varphi ^{*})=\nu (\varphi ).$

We will also use a minor variation on the previous result.

\vglue12pt

\proclaim{Lemma}
\label{P3.2a} Suppose that $L$ and $W$ are operator spaces with $L$ 
finite\/{\rm -}\/dimensional{\rm ,} and let $\iota _{W}:W\rightarrow W^{**}$ denote the canonical
complete isometry{\rm .} Then for any mapping $\varphi :L\rightarrow W,$ we have that
$$
\nu (\iota _{W}\circ \varphi )=\nu (\varphi ).
$$
\endproclaim

\demo{Proof}%
We have a commutative diagram
\begin{eqnarray}
&&\\
\noalign{\vskip-16pt}
&&
\begin{array}{ccc}
L^{*}\hat{\otimes}W & \rightarrow & L^{*}\hat{\otimes}W^{**} \\
\downarrow &  & \downarrow \\
{\cal N}(L,W) & \rightarrow & {\cal N}(L,W^{**})
\end{array}
\nonumber
\end{eqnarray}
where from Lemma~3.1, the top row is a completely isometric
injection, and since $L$ is finite-dimensional, the columns are (complete)
isometries. It follows that the bottom row is a (completely) isometric
injection.%
\enddemo %
 
\section{Completely integral mappings}

As in Banach space theory, the completely nuclear norm is not local. By this
we mean that if we are given operator spaces $V$ and $W,$ and a linear
mapping $\varphi :V\rightarrow W$ such that $\nu (\varphi |_{F})\le 1$ for
all finite-dimensional subspaces $F\subseteq V,$ it need not follow that $%
\varphi $ is completely nuclear. As we will see, this naturally leads to the
more general class of completely integral mappings.

We recall from \cite{ER3} that a linear mapping $\varphi :V\rightarrow W$ is
{\it completely integral} with {\it completely integral norm} $\iota
(\varphi )\le 1$ if $\varphi $ is in the point-norm closure of the set of
finite rank mappings $\psi :V\to W$ such that $\nu (\psi )<1,$ or using a
standard convexity argument (see \cite[Prop.~3.2]{ER3}), $\varphi $ is
in the point-weak closure of that set. We let ${\cal  I}(V,W)$ denote the
linear space of all completely integral mappings from $V$ into $W$ with the
norm $\iota$ and, as usual, we write $\iota (\varphi )=\infty $ if $\varphi
:V\rightarrow W$ is not completely integral. It is clear that we have that $
\iota (\varphi )\leq \nu (\varphi )$ for any linear mapping $\varphi
:V\rightarrow W,$ and thus we have a natural contraction
$$
{\cal  N}(V,W)\rightarrow {\cal  I}(V,W).
$$

\proclaim{Lemma}
\label{P4.2}If $V$ is finite-dimensional{\rm ,} we have the isometry
$$
{\cal  N}(V,W)={\cal  I}(V,W).
$$
\endproclaim

\demo{Proof}%
We must show that if $\varphi :V\rightarrow W$ satisfies $\iota (\varphi
)\le 1,$ then $\nu (\varphi )\le 1.$ Therefore, let us suppose $\varphi $ is
a point-norm limit of mappings $\varphi _{\lambda }\in {\cal  N}(V,W)$
with $\nu (\varphi _{\lambda })<1.$ We fix a basis $x_{1},\ldots ,x_{n}$ for
$V$ and we let $f_{1},\ldots ,f_{n}$ be the corresponding dual basis for $%
V^{*}$; i.e., we define
$$
f_{i}(\sum_{j=1}^{n}c_{j}x_{j})=c_{i}.
$$
Using the algebraic identification ${\cal  CB}(V,W)=V^{*}\otimes W,$ we
have that
$$
\varphi _{\lambda }=\sum_{i=1}^{n}f_{i}\otimes y_{i}^{\lambda },\,\,
$$
and
$$
\varphi =\sum_{i=1}^{n}f_{i}\otimes y_{i},
$$
where $y_{i}^{\lambda }=\varphi _{\lambda }(x_{i})$ and $y_{i}=\varphi
(x_{i})\in W.$ Since $\varphi _{\lambda }$ converges to $\varphi $ in the
point-norm topology, it follows that
$$
\left\| y_{i}^{\lambda }-y_{i}\right\| =\left\| \varphi _{\lambda
}(x_{i})-\varphi (x_{i})\right\| \rightarrow 0.
$$
The operator projective tensor norm is a cross norm in the Banach space
sense, and thus
\begin{eqnarray*}
\nu (\varphi -\varphi _{\lambda }) &\le &\left\| \sum_{i=1}^{n}f_{i}\otimes
(y_{i}^{\lambda }-y_{i})\right\| _{V^{*}\hat{\otimes}W} \\
&\le &\sum_{i=1}^{n}\left\| f_{i}\right\| \left\| y_{i}^{\lambda
}-y_{i}\right\| \rightarrow 0.
\end{eqnarray*}
Since $\nu (\varphi _{\lambda })<1$ and
$$
\nu (\varphi )\le \nu (\varphi -\varphi _{\lambda })+\nu (\varphi _{\lambda
}),
$$
we conclude that $\nu (\varphi )\le 1.$
\enddemo %

Given operator spaces $V$ and $W,$ the pairing (\ref{traceduality}) is given
by
\begin{equation}
\langle \cdot ,\cdot \rangle :{\cal  {\cal  CB}}(V,W)\times (V\otimes
W^{*})\rightarrow \Bbb{C}:(\varphi ,(v\otimes g))\rightarrow \langle \varphi
(v),g\rangle ,\hskip.5in  \label{P4.1a}
\end{equation}
and it thus determines a linear mapping
\begin{equation}
\Psi :{\cal  {\cal  CB}}(V,W)\hookrightarrow (V\otimes W^{*})^{*},
\label{P4.1b}
\end{equation}
where we let $(V\otimes W^{*})^{*}$ denote the space of linear functionals $f$ for
which $f(v\otimes g)$ is norm-continuous in each variable. In particular, since
$$
(V\hat{\otimes}W^{*})^{*}={\cal  CB}(V,W^{**}),
$$
$\Psi $ induces the (completely) isometric injection
\begin{equation}
\Psi :{\cal  {\cal  CB}}(V,W)\hookrightarrow (V\hat{\otimes}W^{*})^{*}
\label{compisom}
\end{equation}
corresponding to the usual inclusion mapping ${\cal  C}{\cal  B}%
(V,W)\subseteq {\cal  CB}(V,W^{**}).$

Modifying  (3.12) in \cite{ER3}, we have a natural commutative diagram
$$
\begin{array}{ccccc}
V^{*}\hat{\otimes}W & \overset{\Psi _{0}}{\longrightarrow } & (V\check{%
\otimes}W^{*})^{*} & \overset{\Phi _{1}^{*}}{\longrightarrow } & (V\hat{%
\otimes}W^{*})^{*} \\
\downarrow {\scriptstyle \Phi} &  &  &  & \uparrow {\scriptstyle \Psi} \\
{\cal  N}(V,W) &  & \hookrightarrow &  & {\cal  CB}(V,W)
\end{array}
$$
where $\Phi $ and $\Phi _{1}:V\hat{\otimes}W^{*}\rightarrow V\check{\otimes}%
W^{*}$ are the canonical (complete) contractions, and the (complete)
contraction $\Psi _{0}$ is determined by the relation
$$
\Psi _{0}(f\otimes w)(v\otimes g)=f(v)g(w).
$$
Since $\Phi _{1}$ has dense range, $\Phi _{1}^{*}$ is one-to-one. It follows
that $\ker \Phi \subseteq \ker \Psi _{0},$ and thus $\Psi _{0}$ determines a
complete contraction $\Psi :{\cal  N}(V,W)\rightarrow (V\check{\otimes}%
W^{*})^{*},$ which is just the restriction of (\ref{compisom}) to ${\cal  N%
}(V,W)$.

\proclaim{Theorem}
\label{compint}
Suppose that $V$and $W$ are operator spaces and that $\varphi
:V\rightarrow W$ is a completely bounded linear mapping{\rm .} Then the following
are equivalent\/{\rm :}\/

\begin{itemize}
\item[{\rm (a)}]  $\iota (\varphi )\le 1,$

\item[{\rm (b$_{{1}}$)}]  $\nu (\varphi _{|E})\le 1 $ for all finite\/{\rm -}\/dimensional subspaces $E\subseteq V,$

\item[{\rm (b$_{{2}}$)}]  $\nu (\varphi \circ \psi )\le 1$ for all complete
contractions $\psi :E\rightarrow V,$ with $E$ finite-dimensional{\rm ,}

\item[{\rm (c$_{{1}}$)}]  $\left\| {\rm id}\otimes \varphi :F\check{\otimes}%
V\rightarrow F\hat{\otimes}W\right\| \le 1$ for  finite-dimensional
operator spaces~$F,$

\item[{\rm (c$_{{2}}$)}]  $\left\| {\rm id}\otimes \varphi :F\check{\otimes}%
V\rightarrow F\hat{\otimes}W\right\| \le 1$ for arbitrary operator spaces $F,
$

\item[{\rm (d)}]  $\left\| \Psi (\varphi ):V\otimes _{\vee }W^{*}\rightarrow \Bbb{C%
}\right\| \le 1.$
\end{itemize}
\endproclaim   

\demo{Proof}%
 (a)$\Leftrightarrow $(b$_{{1}}$). From \cite[Prop.~2.1]{ER5}, we see  that $\iota (\varphi )\le 1$ if and only if
$\iota (\varphi |_{E})\le 1$ for all finite-dimensional subspaces $E\subseteq V.$ Thus the equivalence
follows from Lemma \ref{P4.2}.
\smallbreak
(b$_{{1}}$)$\Leftrightarrow $(b$_{{2}}$). Given $\psi
:E\rightarrow V,$ we have that
$$
\nu (\varphi \circ \psi )\le \nu (\varphi _{|\psi (E)})\left\| \psi \right\|
_{cb}
$$
and the equivalence is immediate.
\smallbreak
(b$_{{2}}$)$\Leftrightarrow $(c$_{{1}}$). Given a finite-dimensional operator space $F$ and letting\break $E=F^{*},$ we
may identify $F%
\check{\otimes}V$ with ${\cal  C}{\cal  B}(E,V).$ This equivalence is
immediate from the commutative diagram
\begin{eqnarray}
&& \label{F2.1a}\\
\noalign{\vskip-16pt}
&&
\begin{array}{ccc}
F\check{\otimes}V & \overset{{\rm id}\otimes \varphi }{\longrightarrow } & F\hat{%
\otimes}W \\
|| &  & || \\
{\cal  CB}(E,V) & \overset{\psi \mapsto \varphi \circ \psi }{%
\longrightarrow } & {\cal  N}(E,W)\ .
\end{array} \nonumber
\end{eqnarray}
\smallbreak
(c$_{{1}}$)$\Leftrightarrow $(c$_{{2}}$). Given an arbitrary
operator space $F$ and an element\break $u\in F\otimes _{\vee }V,$ we have that $%
u\in F_{0}\otimes _{\vee }V$ for some finite-dimensional subspace $F_{0}$ of
$V,$ and
$$
\left\| u\right\| _{F\check{\otimes}V}=\left\| u\right\| _{F_{0}\check{%
\otimes}V}.
$$
Since $F_{0}\hat{\otimes}W\rightarrow F\hat{\otimes}W$ is a contraction, it
is evident that (c$_{{1}}$)$\Rightarrow $(c$_{{2}}$), and the
converse is trivial.

(a)$\Leftrightarrow $(d). Given $\varphi \in {\cal  C}{\cal  B}(V,W)$
with $\iota (\varphi )\le 1$ and an element\break $u\in V\otimes _{\vee }W^{*}$
with $\left\| u\right\| _{\vee }\le 1,$ we may assume that $u\in E\otimes
_{\vee }W^{*},$ where $E$ is a finite-dimensional subspace of $V$. It
follows that
$$
\left| \langle \Psi (\varphi ),u\rangle \right| =\left| \langle \Psi
(\varphi _{|E}),u\rangle \right| \le \nu (\varphi _{|E})\le \iota (\varphi
),
$$
and thus $\left\| \Psi (\varphi )\right\| \le 1$.

Let us suppose that $\varphi \in {\cal  CB}(V,W)$ satisfies $\left\| \Psi
(\varphi )\right\| \le 1.$ We have that
$$
V\check{\otimes}W^{*}\cong W^{*}\check{\otimes}V\subseteq {\cal  CB}%
(W,V^{**})=(V^{*}\hat{\otimes}W)^{*},
$$
and thus $\Psi (\varphi )$ has a contractive extension $\bar{\Psi}_{\varphi
}\in (V^{*}\hat{\otimes}W)^{**}.$ From the\break bipolar theorem, we may choose a
net of elements $u_{\lambda }\in V^{*}\hat{\otimes}W$ such that\break $\left\|
u_{\lambda }\right\| _{V^{*}\hat{\otimes}W}<1$ and $u_{\lambda }$ converges
to $\bar{\Psi}_{\varphi }$ in the weak$^{*}$ topology. Let $\varphi
_{\lambda }=\Phi (u_{\lambda })\in {\cal  N}(V,W)$. Then $\nu (\varphi
_{\lambda })<1$ and
$$
\varphi _{\lambda }(v)(g)=\langle u_{\lambda },v\otimes g\rangle \rightarrow
\langle \bar{\Psi}_{\varphi },v\otimes g\rangle =\langle \varphi ,v\otimes
g\rangle =\langle \varphi (v),g\rangle
$$
for all $v\in V$ and $g\in W^{*}$. Therefore, $\varphi _{\lambda }$
converges to $\varphi $ in the point-weak topology, and $\iota (\varphi )\le
1$.
\enddemo %

\vglue6pt

\proclaim{{C}orollary}
Given operator spaces $V$ and $W$ and a linear mapping $\varphi
:V\rightarrow W,$ we have
\begin{eqnarray}
\noalign{\vskip4pt}
\iota (\varphi ) &=&\sup \left\{ \nu (\varphi _{|E}):E\subseteq V\hbox{ 
finite-dimensional}\right\}    \label{F0} \\ \noalign{\vskip4pt}
&=&\sup \left\{ \left\| \hbox{\rm id}\otimes \varphi :F\check{\otimes}V\rightarrow F%
\hat{\otimes}W\right\| :F\hbox{  finite-dimensional}\right\}   \nonumber \\ \noalign{\vskip4pt}
&=&\sup \left\{ \left\| {\rm id}\otimes \varphi :F\check{\otimes}V\rightarrow F%
\hat{\otimes}W\right\| :F\hbox{  arbitrary}\right\} . \nonumber \\
\noalign{\vskip-9pt}
\nonumber
\end{eqnarray}
Furthermore{\rm ,} the mapping
\begin{equation}
\Psi :{\cal  I}(V,W)\hookrightarrow (V\otimes _{\vee }W^{*})^{*}=\left[
{\cal  CB}^{0}(W,V)\right] ^{*}  \label{isom}
\end{equation}
is an isometric injection{\rm .}
\endproclaim

In particular, we see that $\iota $ is {\it local. }If $W$ is finite-dimensional, 
we have from (\ref{isom}) the isometry
\begin{equation}
{\cal  I}(V,W)=(V\check{\otimes}W^{*})^{*}={\cal  CB}(W,V)^{*}.
\label{integdual}
\end{equation}
However, in contrast to the situation for Banach spaces (see (\ref{classical}%
)), we {\it need not have} that the natural mapping
$$
\tilde{\Psi}:{\cal  I}(V,W^{*})\rightarrow (V\check{\otimes}W)^{*}
$$
is isometric. Using the identification ${\cal  N}(V,W)=V^{*}\hat{\otimes}W$%
, we obtain the following result from the discussion of (\ref{F1.1}).

\proclaim{Proposition}
\label{P4.3a} An operator space $V$ is locally reflexive if and only if we
have the isometry
\begin{equation}
{\cal  N}(V,W)={\cal  I}(V,W)  \label{F4.c}
\end{equation}
for all finite-dimensional $W${\rm .}
\endproclaim 

It is shown in \cite{EH} that $C^{*}(\Bbb{F}_{2})$, the full group $C^{*}$%
-algebra of 2-generators, is not locally reflexive, and thus there exists a
finite-dimensional operator space $W$ such that ${\cal  N}(C^{*}(\Bbb{F}%
_{2}),W)\rightarrow {\cal  I}(C^{*}(\Bbb{F}_{2}),W)$ is not isometric.
Denoting the corresponding Banach mapping spaces with the subscript $B,$ we
always have the isometry
$$
{\cal  N}_{B}(V,W)={\cal  I}_{B}(V,W)
$$
for finite-dimensional $W.$

We conclude this section with a factorization which characterizes the
completely integral mappings.

\proclaim{Proposition}
\label{P4.3b} Given operator spaces $V$ and $W$ and a completely bounded
mapping $\varphi :V\to W,$ we have that $\iota (\varphi )\le 1$ if and only
if there exist Hilbert spaces $H$ and $K,$ a contractive functional $\omega
\in B(H\otimes K)^{*}$ and completely contractive maps $r:V\to B(H)${\rm ,} $%
t:W^{*}\to B(K)$ such that for all $v\in V$ and $g\in W^{*},$
\begin{equation}
\langle \varphi (v),g\rangle =\langle \omega ,r(v)\otimes t(g)\rangle .
\label{E1}
\end{equation}
\endproclaim 


\demo{Proof}%
Let us suppose that $\iota (\varphi )\le 1.$ We fix completely isometric
embeddings $r:V\rightarrow B(H)$ and $s:W^{*}\rightarrow B(K)$. From (4.5)
$$
\left\| \Psi (\varphi ):V\otimes _{\vee }W^{*}\rightarrow \Bbb{C}\right\|
\le 1,
$$
hence we may extend $\Psi (\varphi )$ to an element $\omega \in B(H\otimes K)^{*}$
with $\left\| \omega \right\| \le 1.$ It follows that
$$
\langle \varphi (v),g\rangle =\Psi (\varphi )(v\otimes g)=\omega
(r(v)\otimes t(g)).
$$

Conversely given such a factorization with $\left\| \omega \right\| \left\|
r\right\| _{cb}\left\| t\right\| _{cb}\le 1$, we have that for any $u\in
V\otimes W^{*},$%
\begin{eqnarray*}
\noalign{\vskip4pt}
\left| \Psi (\varphi )(u)\right| &=&\left| \langle \omega ,(r\otimes
t)(u)\rangle \right| \\ \noalign{\vskip4pt}
&\le &\left\| \omega \right\| \left\| (r\otimes t)(u)\right\| _{B(H\otimes
K)} \\ \noalign{\vskip4pt}
&\le &\left\| \omega \right\| \left\| r\right\| _{cb}\left\| t\right\|
_{cb}\left\| u\right\| _{V\otimes _{\vee }W^{*}} \\ \noalign{\vskip4pt}
&\le &\left\| u\right\| _{V\otimes _{\vee }W^{*}}.\\
\noalign{\vskip-9pt}
\end{eqnarray*}
Therefore, we have $\iota (\varphi )=\left\| \Psi (\varphi )\right\| \le 1$.
\enddemo %

Given a bounded linear functional $\omega :B(H\otimes K)\rightarrow \Bbb{C}$%
, we define a linear mapping
$$
M(\omega ):B(H)\rightarrow B(K)^{*}
$$
by
$$
M(\omega )(b)(g)=\omega (b\otimes g).
$$

\proclaim{{C}orollary}
\label{P4.3c}
Let us suppose that $V$ and $W$ are operator spaces{\rm ,} and that $%
\varphi :V\rightarrow W$ is a linear mapping{\rm .} We have that $\iota (\varphi
)\le 1$ if and only if there is a commutative diagram
\begin{eqnarray}
&& \label{E2}\\
\noalign{\vskip-18pt}
&&
\begin{array}{ccccc}
B(H) & \overset{M(\omega )}{\longrightarrow } & B(K)^{*} &  &  \\
{\scriptstyle \ r}\uparrow  &  &  & \searrow {\scriptstyle s} &  \\
V & \overset{\varphi }{\longrightarrow } & W & \overset{\iota _{W}}{%
\hookrightarrow } & W^{**}\ ,  
\end{array} \nonumber
\end{eqnarray}
where $\omega \in B(H\otimes K)^{*}\ $ satisfies $\left\| \omega \right\|
\le 1,$ $r$ and $s$ are complete contractions{\rm ,} $\iota _{W}:W\rightarrow
W^{**}$ is the canonical embedding{\rm ,} and $s:B(K)^{*}\rightarrow W^{**}$ is
weak$^{*}$ continuous{\rm .}
\endproclaim

\demo{Proof}%
Letting $s=t^{*},$ this is immediate from Proposition \ref{P4.3b}.%
\enddemo %

\section{Exactly integral mappings}

Weakening the characterization in Corollary \ref{P4.3c}, we say that a
linear mapping $\varphi :V\rightarrow W$ is {\it exactly integral} if it
has a factorization (\ref{E2}), where $r$ and $s$ are completely bounded and
$\omega \in B(H\otimes K)^{*}$, {\it but we do not assume that }$s${\it \
is weak}$^{*}${\it \ continuous}. We define the corresponding {\it exactly
integral norm}
$$
\iota^{\rm ex}(\varphi )={\rm inf}\{\left\| r\right\| _{cb}\left\| \omega
\right\| \left\| s\right\| _{cb}\}
$$
where the infimum is taken over all such factorizations. It is trivial that
if $\varphi :V\rightarrow W$ is completely integral, then $\varphi $ is
exactly integral and $\iota^{\rm ex}(\varphi )\le \iota (\varphi ).$ The fact
that $\iota ^{\rm ex}$ is a norm follows from Theorem \ref{P4.4a}.

\proclaim{Lemma}
\label{P4.4} Let us suppose that $V$ and $W$ are operator spaces{\rm .} If\break $%
\varphi :V\rightarrow W$ is completely  integral{\rm ,} then $\varphi ^{*}:W^{*}\rightarrow
V^{*}$ is exactly integral with $\iota ^{\rm ex}(\varphi ^{*})\break\le \iota (\varphi).$
\endproclaim

\demo{Proof}%
We may use (\ref{E2}) to construct a commutative diagram
$$
\begin{array}{ccccc}
B(K) & \overset{M(\tilde{\omega})}{\longrightarrow } & B(H)^{*} &  &  \\
{\scriptstyle s_{*}}\uparrow &  &  & \searrow {\scriptstyle \iota
_{V^{*}}\circ r^{*}} &  \\
W^{*} & \overset{\varphi ^{*}}{\longrightarrow } & V^{*} & \overset{\iota
_{V^{*}}}{\hookrightarrow } & V^{***}\ ,
\end{array}
$$
where $s=(s_{*})^{*},$ and $\tilde{\omega}:B(K\otimes H)\rightarrow \Bbb{C}$
is the obvious ``flip'' of $\omega .$%
\enddemo %

It will be noted that in the above proof, $\iota _{V^{*}}\circ r^* $ is
generally not weak$^{*}$ continuous, and thus we cannot conclude that $\iota
(\varphi ^{*})\le \iota (\varphi ).$

\proclaim{Lemma}
\label{P4.5} If $A$ is a $C^{*}$\/{\rm -}\/algebra and $V$ is an arbitrary operator
space{\rm ,} then we have the isometric identification
$$
{\cal  I}^{\rm ex}(V,A)={\cal  I}(V,A).
$$
\endproclaim

\demo{Proof}%
Let us assume that $\iota ^{\rm ex}(\varphi )\le 1$. Then we can find a
factorization
$$
\begin{array}{ccccc}
B(H) & \overset{M(\omega )}{\longrightarrow } & B(K)^{*} &  &  \\
{\scriptstyle r}\uparrow &  &  & \searrow {\scriptstyle \ s} &  \\
V & \overset{\varphi }{\longrightarrow } & A & \overset{\iota _{A}}{%
\hookrightarrow } & A^{**},
\end{array}
$$
where $r,s$ are complete contractions and $\omega $ is of norm one. From
Theorem \ref{P2.1} we may approximate $s$ in the point-weak$^{*}$ topology
by a net of weak$^{*}$ continuous mappings $s_{\lambda }:B(K)^{*}\rightarrow
A$ with $\left\| s_{\lambda }\right\| _{cb}\le 1$. Fixing $\lambda ,$ and
letting $\varphi _{\lambda }=\iota _{A}s_{\lambda }M(\omega )r,$ we have the
commutative diagram
$$
\begin{array}{ccccc}
B(H) & \overset{M(\omega )}{\longrightarrow } & B(K)^{*} &  &  \\
{\scriptstyle \ r}\uparrow &  & \downarrow {\scriptstyle s_{\lambda }} &
\searrow \scriptstyle\iota _{A}\circ {s}_{\lambda } &  \\
V & \overset{\varphi _{\lambda }}{\longrightarrow } & A & \overset{\iota _{A}%
}{\hookrightarrow } & A^{**},
\end{array}
$$
where $\iota _{A}\circ s_{\lambda }:B(K)^{*}\rightarrow A^{**}$ is a weak$%
^{*}$-continuous complete contraction. It follows from Corollary \ref{P4.3c}
that $\iota (\varphi _{\lambda })\le 1.$ Since each $s_{\lambda }$ and $%
\varphi $ have range in $A,$ $\varphi _{\lambda }$ converges to $\varphi $
in the point-weak topology. Thus we have  from the definition of the completely
integral norm that $\iota (\varphi )\le 1.$
\enddemo %

Although the definition of the exactly integral mappings might seem
contrived, such mappings play a natural and important role in operator space
theory. In order to substantiate this claim, we will provide several
alternative characterizations. This material will not be needed in the
subsequent sections.

The following is well-known:

\proclaim{Lemma}
\label{L1} Suppose that $E$ is a matrix space {\rm (}\/see \S {\rm 1).} Then for any
operator space $W,$ we have the complete isometry
\begin{equation}
(E\check{\otimes}W)^{*}\cong E^{*}\hat{\otimes}W^{*}.  \label{matrix}
\end{equation}
\endproclaim

\demo{Proof}%
Let us suppose that $E\subseteq M_{n},$ and let $\rho :M_{n}^{*}\rightarrow
E^{*}$ be the restriction mapping. We have that $E\check{\otimes}W\subseteq
M_{n}\check{\otimes}W,$ and this determines the restriction mapping $\bar{%
\rho}$ in the  commutative diagram
$$
\begin{array}{ccc}
T_{n}\hat{\otimes}W^{*} & \cong & (M_{n}\check{\otimes}W)^{*} \\
{\scriptstyle \rho \otimes id}\downarrow &  & {\scriptstyle \bar{\rho}}%
\downarrow \\
E^{*}\hat{\otimes}W^{*} &\rightarrow & (E\check{\otimes}W)^{*}\ .
\end{array}
$$
From the general theory, it follows that the top row is completely isometric,
and the first column is a complete quotient mapping. On the other hand,
since $E\check{\otimes}W\rightarrow M_{n}\check{\otimes}W$ is a complete
isometry, the second column is a complete quotient mapping. It follows that
the bottom row is a complete isometry.%
\enddemo %

By contrast to the situation for Banach spaces, if $E$ is a 
finite-dimensional operator space, (\ref{matrix}) need not hold in general. This is
related to the fact that $E$ need not be exact, i.e., approximable in the
Pisier-Banach-Mazur sense by matrix spaces (see \cite{Pi0}). Nevertheless,
it can be approximated in an asymptotic sense. We may identify $E$ with a
subspace of $M_{\infty }.$ For each $n\in \Bbb{N},$ we let
$$
P_{n}:x\in M_{\infty }\rightarrow x^{(n)}\in M_{n}
$$
be the usual truncation mapping. Restricting both the domain and the range,
we let
\begin{equation}
\tau _{n}=\tau _{n}^{E}=P_{n|E}:E\rightarrow P_{n}(E).  \label{taudef}
\end{equation}

\proclaim{Lemma}
\label{L2} Given a finite-dimensional subspace $E$ of $M_{\infty },$ an
integer $k>0,$ and $0<\varepsilon <1,$ there exists an $n\in \Bbb{N}$ such that $%
\tau _{n}$ is invertible and
$$
\left\| (\tau _{n})_{k}^{-1}\right\| \leq 1+\varepsilon .
$$
\endproclaim

\demo{Proof}%
Let us fix elements $x_{i}$ which are $\varepsilon /2$-dense in the unit sphere
of $E.$ Since
$$
\lim_{n\rightarrow \infty }\left\| P_{n}(x_{i})\right\| =\left\|
x_{i}\right\| ,
$$
we may choose an $n$ such that
$$
\left\| P_{n}(x_{i})\right\| \geq 1-\varepsilon /2
$$
for all $i.$ If $x\in E$ and $\left\| x\right\| =1,$ we may find an $i$ such
that $\left\| x-x_{i}\right\| <\varepsilon /2.$ It follows that
$$
\left\| P_{n}(x)\right\| \geq \left\| P_{n}(x_{i})\right\| -\left\|
P_{n}(x_{i})-P_{n}(x)\right\| \geq 1-\varepsilon,
$$
and thus
$$
\left\| \tau _{n}(x)\right\| \geq (1-\varepsilon )\left\| x\right\| .
$$
It follows that $\tau _{n}$ is one-to-one, and $\left\| \tau
_{n}^{-1}\right\| \leq (1-\varepsilon )^{-1}.$ We may apply this argument to
the mappings
$$
(P_{n})_{k}:M_{k}(E)\rightarrow M_{k}(P_{n}(E)),
$$
and the result follows.%
\enddemo %

We assume that readers are familiar with ultraproducts of Banach spaces and
operator spaces (see [12], [15], [28], [29], and [35]).

Groh has proved that an ultrapower of von Neumann algebraic preduals is
again the predual of a von Neumann algebra (see \cite{UG} and \cite{Ro} ---
we are indebted to Ward Henson for bringing these papers to our attention).
In order to make our discussion more accessible, we repeat his argument for
the von Neumann algebra $M_{\infty }.$ Given an index set $I,$ and a free
ultrafilter ${\cal  U}$ on $I,$ we let $\prod_{{\cal  U}}T_{\infty }$
denote the operator space ultrapower of $T_{\infty }$. We have a natural
completely isometric injection
\begin{eqnarray}
\noalign{\vskip4pt}
&&
\theta :\prod_{{\cal  U}}T_{\infty }\rightarrow \ell ^{\infty
}(I,M_{\infty })^{*}  \label{embed}\\
\noalign{\vskip-9pt}\nonumber
\end{eqnarray}
defined by the pairing
\begin{eqnarray*}
\noalign{\vskip4pt}
&&
\left\langle \theta ([\omega _{\alpha }]),(y_{\alpha })\right\rangle =\lim_{%
{\cal  U}}\left\langle \omega _{\alpha },y_{\alpha }\right\rangle
\\
\noalign{\vskip-9pt}
\end{eqnarray*}
(see, e.g., \cite{ER5}). We may regard $\ell ^{\infty }(I,M_{\infty })^{*}$
as a bimodule over $\ell ^{\infty }(I,M_{\infty })$ or over $\ell ^{\infty
}(I,M_{\infty })^{**}$ in the usual manner. The subspace $T=\theta (\prod_{%
{\cal  U}}T_{\infty })$ is a norm closed two-sided $\ell ^{\infty
}(I,M_{\infty })$ submodule since if we are given $f=[(f_{\alpha })]\break\in
\prod_{{\cal  U}}T_{\infty }$ and $x=(x_{\alpha })\in \ell ^{\infty
}(I,M_{\infty }),$ we have that $(x_{a}f_{\alpha })\in \ell ^{\infty
}(I,T_{\infty })$ and thus
\begin{eqnarray*}
&&
xf=\theta [(x_{a}f_{\alpha })]\in T,
\\
\noalign{\vskip-9pt}
\end{eqnarray*}
and the same argument shows that $fx\in T$. We conclude (see \cite[Chap.~III, Th.~2.7]{Ta}) that the annihilator of $T$
is a weak$^{*}$ closed two-sided ideal in the von Neumann algebra $\ell ^{\infty }(I,M_{\infty
})^{**}$, and in particular, there is a central projection $e\in \ell
^{\infty }(I,M_{\infty })^{**}$ for which
\begin{eqnarray}
\noalign{\vskip4pt}
&&
T=\ell ^{\infty }(I,M_{\infty })^{*}e=\left[ \ell ^{\infty }(I,M_{\infty
})^{**}e\right] _{*}.  \label{vonNeum}
\\
\noalign{\vskip-9pt}\nonumber
\end{eqnarray}

It is useful to compare the following theorem with Theorem \ref{compint}.
Significant portions of this result may be found in \cite{J}, where a rather
different approach is used. Condition (d) is related to Pisier's
factorization theorem for completely 1-summing mappings \cite{Pi1}.

\proclaim{Theorem}
\label{P4.4a} Given operator spaces $V$ and $W$ and a completely bounded
mapping $\varphi :V\rightarrow W${\rm ,} the following are equivalent\/{\rm :}\/

\begin{itemize}
\item[{\rm (a)}]  $\iota ^{\rm ex}(\varphi )\le 1,$

\item[{\rm (b)}]  $\nu (\varphi \circ \psi )\le 1$ for all complete contractions $%
\psi :E\rightarrow V$ with $E$ a matrix space{\rm ,}

\item[{\rm (c)}]  $\left\| {\rm id}\otimes \varphi :E^{*}\check{\otimes}V\rightarrow
E^{*}\hat{\otimes}W\right\| \le 1$ for all matrix spaces $E,$

\item[{\rm (d)}]  There exists an infinite index set $I,$ a free ultrafilter $%
{\cal  U}$ on $I,$ and a commutative diagram
\begin{eqnarray}
&&\label{ultradiag}
\\
\noalign{\vskip-16pt}
&&
\begin{array}{ccccc}
\ell ^{\infty }(I,M_{\infty }) & \overset{{\cal  M}}{\longrightarrow } &
\prod_{{\cal  U}}T_{\infty } &  &  \\
{\scriptstyle r}\uparrow  &  &  & \searrow {\scriptstyle s} &  \\
V & \overset{\varphi }{\longrightarrow } & W & \overset{\iota _{W}}{%
\hookrightarrow } & W^{**}\ ,  
\end{array}
\nonumber
\end{eqnarray}
where $r$ and $s$ are complete contractions{\rm ,} and ${\cal  M}=[M(a_{\alpha
},b_{\alpha })]$ is determined by the multiplication operators $M(a_{\alpha
},b_{\alpha }):M_{\infty }\rightarrow T_{\infty }$ with $\left\| a_{\alpha
}\right\| _{2},\left\| b_{\alpha }\right\| _{2}< 1.$
\end{itemize}
\endproclaim 

\demo{Proof}%
(a)$\Rightarrow $(b). Let us suppose that we have a factorization (\ref{E2})
for the mapping $\varphi ,$ with $\omega $ contractive, and $r$ and $s$
completely contractive. Given a matrix space $E$ and a complete contraction $%
\psi :E\rightarrow V,$ we have from Lemma 5.3 that  $\omega ((r\circ \psi )\otimes id)$ is a
strictly
contractive element of
$$
(E\check{\otimes}B(K))^{*}=E^{*}\hat{\otimes}B(K)^{*}.
$$
The corresponding element of ${\cal N}(E,B(K)^{*})$ is $M(\omega )\circ r\circ
\psi ,$ since if $x\in E,$ and $b\in B(K),$%
$$
M(\omega )(r(\psi (x))(b)=\omega (r(\psi (x))\otimes b)=\omega ((r\circ \psi
)\otimes id)(x\otimes b)).
$$
Thus using Lemma \ref{P3.2a} and the factorization (\ref{E2}),
$$
\nu (\varphi \circ \psi )=\nu (\iota _{W}\circ \varphi \circ \psi )=\nu
(s\circ M(\omega )\circ r\circ \psi )\le \left\| s\right\| _{cb}\nu
(M(\omega )\circ r\circ \psi )\le 1.
$$
\medbreak
(b)$\Leftrightarrow $(c) is immediate from the commutative diagram (\ref
{F2.1a}).
\medbreak
(c)$\Rightarrow $(d) We let $I$ be the index set of all triples $\alpha
=(E,F,k)$, where $E\subseteq V$ is finite-dimensional, $F\subseteq W$ is
finite-codimensional, and $k\in \Bbb{N}$. Given such a triple $\alpha ,$ we
shall also use the notation $E=E_{\alpha },F=F_{\alpha },$ and $k=k_{\alpha
}.$ We write $\iota _{\alpha }:E_{\alpha }\hookrightarrow V$ and $\pi
_{\alpha }:W\rightarrow W/F_{\alpha }$ for the inclusion and quotient
mappings. We define a partial order on $I$ by $\alpha \preceq \alpha
^{\prime }=(E^{\prime },F^{\prime },k^{\prime })$ if $E\subseteq E^{\prime }$%
, $F^{\prime }\subseteq F$, $\ $and $k\le k^{\prime }$. For each $\alpha \in
I$, we let $I_{\alpha }=\{\alpha ^{\prime }\in I:\alpha \preceq \alpha
^{\prime }\}$, we write ${\cal  F}_{\preceq }$ for the filter generated by
these $I_{\alpha }$'s and we fix a free ultrafilter ${\cal  U}$ on $I$\break
containing~${\cal  F}_{\preceq }.$

For each $\alpha =(E,F,k)\in I$, $W/F$ is a finite-dimensional operator
space, and thus we may identify it with a finite-dimensional subspace $%
G=G_{\alpha }$ of $M_{\infty}$, and for each $n\in \Bbb{N},$ we let $\tau _{n}^{G_{\alpha }}=
P_{n|G_\alpha}$ (see (\ref{taudef})). From
Lemma \ref {L2}, we may choose an integer $n(\alpha )\in \Bbb{N}$ with $\tau _{n(\alpha
)}^{G_{\alpha }}$ invertible and with
\begin{equation}
\left\| (\tau _{n(\alpha )}^{G_{\alpha }})_{k(\alpha )}^{-1}\right\|< 1+%
\frac{1}{k(\alpha )}.  \label{inverse}
\end{equation}
We can choose a constant $0<c_\alpha<1$ so that  $\tau _{\alpha }=c_\alpha \tau^{G_\alpha}_{n(\alpha)}$ 
also satisfies $\Vert (\tau_\alpha)^{-1}_{k(\alpha)}\Vert <1+{1\over k(\alpha)}$.
We have that

\begin{eqnarray*}
\noalign{\vskip-36pt}\\
&&\varphi _{\alpha }=\tau _{\alpha }\circ \pi _{\alpha }\circ \varphi \circ
\iota _{\alpha }
\end{eqnarray*}
is a linear mapping from $E _{\alpha }$ onto a matrix space $N_{\alpha
}\subseteq M_{n(\alpha )}$.

From (c) we see  that
$$
\left\| {\rm id}\otimes \varphi :N_{\alpha }^{*}\check{\otimes}V\rightarrow
N_{\alpha }^{*}\hat{\otimes}W\right\| \le 1,
$$
and thus
$$
\left\| {\rm id}\otimes \varphi _{\alpha }:N_{\alpha }^{*}\check{\otimes}E_{\alpha
}\rightarrow N_{\alpha }^{*}\hat{\otimes}N_{\alpha }\right\| < 1.
$$
It follows that if we are given an element $\psi \in N_{\alpha }^{*}\check{%
\otimes}E_{\alpha },$ we have from (\ref{traceduality2}) and (\ref
{traceduality3}) that
$$
|\langle \varphi _{\alpha },\psi \rangle |=\left| \hbox{{\rm trace}}%
\,({\rm id}\otimes \varphi _{\alpha })(\psi )\right| \le \left\| ({\rm id}\otimes
\varphi _{\alpha })(\psi )\right\| < \left\| \psi \right\| _{N_{\alpha
}^{*}\check{\otimes}E_{\alpha }}.
$$
We conclude that $\varphi _{\alpha }$ is a strictly contractive element in
$$
(N_{\alpha }^{*}\check{\otimes}E_{\alpha })^{*}=(N_{\alpha }\hat{\otimes}%
E_{\alpha }^{*})^{**}=N_{\alpha }\hat{\otimes}E_{\alpha }^{*}={\cal  N}%
(E_{\alpha },N_{\alpha })
$$
(the second dual of a finite-dimensional operator space coincides with
itself), and  from (3.3)  we have a commutative diagram
$$
\begin{array}{ccc}
M_{\infty } & \overset{M(a_{\alpha },b_{\alpha })}{\longrightarrow } &
T_{\infty } \\
{\scriptstyle r}_{\alpha }\uparrow &  & \downarrow {\scriptstyle s}_{\alpha }
\\
E_{\alpha } & \overset{\varphi _{\alpha }}{\longrightarrow } & N_{\alpha },
\end{array}
$$
where $r_{\alpha }$ and $s_{\alpha }$ are complete contractions, and $%
\left\| a_{\alpha }\right\| _{2},\left\| b_{\alpha }\right\| _{2}\le 1$.
Letting $\tilde{r}_{\alpha }:V\rightarrow M_{\infty }$ be a completely
contractive extension of $r_{\alpha }$ to $V$, we obtain the following
commutative diagram
$$
\begin{array}{ccccccccc}
&  & M_{\infty } &  & \overset{M(a_{\alpha },b_{\alpha })}{\longrightarrow }
&  & T_{\infty } & \overset{\tau _{\alpha }^{-1}\circ s_{\alpha }}{%
\longrightarrow } & W/F_{\alpha } \\
& {\scriptstyle \tilde{r}}_{\alpha }\nearrow & {\scriptstyle r}_{\alpha
}\uparrow &  &  &  &  & \searrow {\scriptstyle s}_{\alpha } & {\uparrow }{%
\scriptstyle \tau }_{\alpha }^{-1} \\
V & \overset{\iota _{\alpha }}{\hookleftarrow } & E_{\alpha } & \overset{%
\varphi \circ \iota _{\alpha }}{\longrightarrow } & W & \overset{\pi
_{\alpha }}{\longrightarrow } & W/F_{\alpha } & \overset{\tau _{\alpha }}{%
\longrightarrow } & N_{\alpha }\ .
\end{array}
$$

We let $\tilde{r}=(\tilde{r}_{\alpha }):V\rightarrow \ell ^{\infty
}(I,M_{\infty })$, and ${\cal  M}=[M(a_{\alpha },b_{\alpha })]$. The
mappings $s_{\alpha }$ and $\tau _{\alpha }^{-1}$ determine corresponding
ultraproduct mappings
$$
\lbrack \tau _{\alpha }^{-1}\circ s_{\alpha }]:\prod_{{\cal  U}}T_{\infty }%
\overset{[s_{\alpha }]}{\longrightarrow }\prod_{{\cal  U}}N_{\alpha }%
\overset{[\tau _{\alpha }^{-1}]}{\longrightarrow }\prod_{{\cal  U}%
}W/F_{\alpha },
$$
where $[s_{\alpha }]$ is a complete contraction, and from (\ref{inverse}),
it is evident that the same is true for $[\tau _{\alpha }^{-1}]$. Finally,
given $[w_{\alpha }+F_{\alpha }]\in \prod_{{\cal  U}}W/F_{\alpha }$, we
may assume that $(w_{\alpha })$ is a uniformly bounded net of   elements in $W$. Then
the weak$^{*}$ ultralimit $\lim_{{\cal  U}}w_{\alpha }$ exists in $W^{**}$
(the unit ball is weak$^{*}$ compact), and we may define $\tilde{s}%
([w_{\alpha }+F_{\alpha }])=\lim_{{\cal  U}}w_{\alpha }$. It is easy to
see that $\tilde{s}:\prod_{{\cal  U}}W/F_{\alpha }\to W^{**}$ is a
well-defined complete contraction such that the diagram
$$
\begin{array}{ccc}
& \prod_{{\cal  U}}W/F_{\alpha } &  \\
{\scriptstyle [\pi }_{\alpha }{\scriptstyle ]}\nearrow &  & \searrow {%
\scriptstyle \tilde{s}} \\
W & \overset{\iota _{W}}{\longrightarrow } & W^{**}
\end{array}
$$
commutes. This gives us the following commutative diagram of complete
contractions
$$
\begin{array}{ccccc}
\ell ^{\infty }(I,M_{\infty }) & \overset{{\cal  M}}{\longrightarrow } &
\prod_{{\cal  U}}T_{\infty } & \overset{{[\tau _{\alpha }^{-1}\circ
s_{\alpha }]}}{\longrightarrow } & \prod_{{\cal  U}}W/F_{\alpha } \\
{\scriptstyle \tilde{r}}\uparrow &  &  & {\scriptstyle [\pi }_{\alpha }{%
\scriptstyle ]}\nearrow & \downarrow {\scriptstyle \tilde{s}} \\
V & \overset{\varphi }{\longrightarrow } & W & \overset{\iota _{W}}{%
\longrightarrow } & W^{**}\ .
\end{array}
$$
Letting $s=\tilde{s}\circ {[\tau _{\alpha }^{-1}\circ s_{\alpha }],}$ we
obtain (\ref{ultradiag}).
\smallbreak
(d)$\Rightarrow $(a) Assuming that we have the factorization (\ref{ultradiag}%
), our task is to construct from it a factorization of the form (\ref{E2}).
From (\ref{vonNeum}), we have a complete isometry
$$
\theta :\prod_{{\cal  U}}T_{\infty }\rightarrow T=\ell ^{\infty
}(I,M_{\infty })^{*}e,
$$
where $e$ is a central projection in $\ell ^{\infty }(I,M_{\infty })^{**},$
and thus we may define a projection $P_{e}$ of $\ell ^{\infty }(I,M_{\infty
})^{*}$ onto $\prod_{{\cal  U}}T_{\infty }$ by letting $P_{e}(f)=fe.$ We
may use this to elaborate (\ref{ultradiag}) in the commutative diagram
$$
\begin{array}{ccccc}
\ell ^{\infty }(I,M_{\infty }) & \overset{{\cal  M}}{\longrightarrow } &
\prod_{{\cal  U}}T_{\infty } & \overset{\theta }{\longrightarrow } & \ell
^{\infty }(I,M_{\infty })^{*} \\
{\scriptstyle \ r}\uparrow &  &  & \searrow {\scriptstyle s} & \downarrow {%
\scriptstyle s\circ P_{e}} \\
V & \overset{\varphi }{\longrightarrow } & W & \overset{\iota _{W}}{%
\hookrightarrow } & W^{**}\ .
\end{array}
$$

Turning to the left side of this diagram, we may assume that $V$ is an
operator subspace of $B(H)$ and let $\iota :V\hookrightarrow B(H)$ be the
inclusion mapping. We have that $r=(r_{\alpha })$ where each $r_{\alpha
}:V\rightarrow M_{\infty }$ is a complete contraction, and using the
Arveson-Wittstock Hahn-Banach theorem, we may extend each $r_{\alpha }\ $to
a complete contraction $\tilde{r}_{\alpha }:B(H)\rightarrow M_{\infty }$.
These determine the complete contraction $\tilde{r}=(\tilde{r}_{\alpha
}):B(H)\rightarrow \ell ^{\infty }(I,M_{\infty })$, and we have the
commutative diagram.
$$
\begin{array}{ccccccc}
B(H) & \overset{\tilde{r}}{\longrightarrow } & \ell ^{\infty }(I,M_{\infty })
& \overset{{\cal  M}}{\longrightarrow } & \prod_{{\cal  U}}T_{\infty } &
\overset{\theta }{\longrightarrow } & \ell ^{\infty }(I,M_{\infty })^{*} \\
{\scriptstyle \iota }\uparrow & {\scriptstyle r}\nearrow &  &  &  & \searrow
{\scriptstyle s} & \downarrow {\scriptstyle s\circ P_{e}} \\
V &  & \overset{\varphi }{\longrightarrow } & W & \overset{\iota _{W}}{%
\hookrightarrow } &  & W^{**}.
\end{array}
$$

For each $\alpha \in I$, we let $\omega _{\alpha }:B(H)\otimes M_{\infty
}\rightarrow \Bbb{C}$ be the linear functional given by
$$
\omega _{\alpha }(x\otimes y)=\left\langle a_{\alpha }{\tilde{r}}_{\alpha
}(x)b_{\alpha },y\right\rangle .
$$
Since $\left\| a_{\alpha }\right\| _{2},\left\| b_{\alpha }\right\| _{2}\le
1 $, it is clear that $\omega _{\alpha }$ is a contractive linear functional
on $B(H)\check{\otimes}M_{\infty }$. Then $[\omega _{\alpha }]$ is a
contractive element in
$$
\prod_{{\cal  U}}(B(H)\ \check{\otimes} \ M_{\infty })^{*}\subseteq \ell
^{\infty }(I,B(H)\ \check{\otimes} \ M_{\infty })^{*}
$$  
where we have used the corresponding identification of the ultraproduct with
a subspace of $\ell ^{\infty }(I,B(H)\ \check{\otimes} \ M_{\infty })^{*}$ (see (%
\ref{embed})). We can identify $B(H)\ \check{\otimes} \ \ell ^{\infty
}(I,M_{\infty })$ with an operator subspace of $\ell ^{\infty }(I,B(H)\check{%
\otimes}M_{\infty })$, and we let $\omega $ be the restriction of $[\omega
_{\alpha }]$ to $B(H)\ \check{\otimes} \ \ell ^{\infty }(I,M_{\infty })$. Then $%
\omega $ is a contractive linear functional on $B(H)\otimes _{\vee }\ell
^{\infty }(I,M_{\infty })$ such that for every $x\in B(H)$ and $(y_{\alpha
})\in \ell ^{\infty }(I,M_{\infty })$,
$$
\omega (x\otimes (y_{\alpha }))=\lim_{{\cal  U}}\omega _{\alpha }(x\otimes
y_{\alpha })=\lim_{{\cal  U}}\left\langle a_{\alpha }{\tilde{r}}_{\alpha
}(x)b_{\alpha },y_{\alpha }\right\rangle =\left\langle \theta \circ {}%
{\cal  M}\circ {\tilde{r}}(x),(y_{\alpha })\right\rangle .
$$
This shows that $M(\omega )=\theta \circ {\cal  M}{}\circ {\tilde{r}}$.

Finally, we let $J$:$\ell ^{\infty }(I,M_{\infty })\hookrightarrow B(K)$ be
an identification of $\ell ^{\infty }(I,M_{\infty })$ with a von Neumann
subalgebra of $B(K)$ for some Hilbert space $K$. Since $\ell ^{\infty
}(I,M_{\infty })$ is injective, there is a completely contractive projection
$\pi $ from $B(K)$ onto $\ell ^{\infty }(I,M_{\infty }).$ Taking adjoints we
have that the composition
$$
\ell ^{\infty }(I,M_{\infty })^{*}\overset{\pi ^{*}}{\longrightarrow }%
B(K)^{*}\overset{J^{*}}{\longrightarrow }\ell ^{\infty }(I,M_{\infty })^{*}
$$
is just the identity mapping, and we obtain the commutative diagram
$$
\begin{array}{ccccc}
B(H) & \overset{M(\omega )}{\longrightarrow } & \ell ^{\infty }(I,M_{\infty
})^{*} & \overset{\pi ^{*}}{\longrightarrow } & B(K)^{*} \\
{\scriptstyle \iota }\uparrow &  &  &  & \downarrow {\scriptstyle s\circ
P_{e}\circ J^{*}} \\
V & \overset{\varphi }{\longrightarrow } & W & \overset{\iota _{W}}{%
\hookrightarrow } & W^{**}\ .
\end{array}
$$
The composition $\omega \circ ({\rm id}\otimes \pi )$ is a contractive functional
on $B(H)\ \check{\otimes} \ B(K);$  thus we may extend it to a contractive
functional $\tilde{\omega}$ on $B(H\otimes K).$ For any $z\in B(K)$ we have
that
$$
M(\tilde{\omega})(x)(z)=\omega (x\otimes \pi (z))=(M(\omega )(x))(\pi
(z))=(\pi ^{*}M(\omega )(x))(z);
$$
 thus $\pi ^{*}\circ M(\omega )=M(\tilde{\omega}).$ We obtain the
commutative diagram
$$
\begin{array}{ccccc}
B(H) & \overset{M(\tilde{\omega})}{\longrightarrow } & B(K)^{*} &  &  \\
{\scriptstyle \iota }\uparrow &  &  & \searrow {\scriptstyle \tau} &  \\
V & \overset{\varphi }{\longrightarrow } & W & \overset{\iota _{W}}{%
\hookrightarrow } & W^{**},
\end{array}
$$
where $\tau =s\circ P_{e}\circ J^{*}$ is a complete contraction, and we
conclude that $\varphi $ is exactly integral with $\iota ^{\rm ex}(\varphi )\le
1 $.%
\enddemo %

Given operator spaces $V$ and $W$, we let ${\cal  I}^{\rm ex}(V,W)$ denote the
space of all exactly integral mappings from $V$ into $W$. It is easy to see
from Theorem \ref{P4.4a} that ${\cal  I}^{\rm ex}(V,W)$ is an operator ideal,
i.e., $\iota ^{\rm ex}$ satisfies (\ref{opid}). It was shown in \cite{J} that $%
\iota ^{\rm ex}$ can also be characterized as a dual operator norm. We have from
(c) of Theorem \ref{P4.4a} and (\ref{matrix}),
$$
\iota ^{\rm ex}(\varphi )=\sup \left\{ \left| \langle ({\rm id}\otimes \varphi
)(u),v\rangle \right| \right\} ,
$$
where the supremum is taken over all $u\in E^{*}\ \check{\otimes} \ V$ and $v\in E%
\ \check{\otimes} \ W^{*}$ and $\left\| u\right\| ,\left\| v\right\| \le 1$ with $%
E$ an arbitrary matrix space. If we let $u$ and $v$ correspond to the
functions $a\in {\cal  C}{\cal  B}(E,V)$ and $b\in {\cal  C}{\cal  B}%
(W,E),$ a simple calculation with elementary matrices leads to the formula
\begin{equation}
\iota ^{\rm ex}(\varphi )=\sup \left\{ {\rm trace\,}(\varphi \circ \psi
):\psi =a\circ b,\left\| a\right\| _{cb},\left\| b\right\| _{cb}\le
1\right\} .  \label{exactint}
\end{equation}

In general given a finite rank mapping $\psi :W\rightarrow V,$ we define
$$
\gamma _{SK}(\psi )=\inf \left\{ \left\| a\right\| _{cb}\left\| b\right\|
_{cb}\right\}
$$
where the supremum is taken over all factorizations
$$
\begin{array}{ccccc}
&  & E &  &  \\
& {\scriptstyle b}\nearrow &  & \searrow {\scriptstyle a} &  \\
W &  & \overset{\psi }{\longrightarrow } &  & V
\end{array}
$$
with $E$ a matrix space. It is easy to see that this determines a norm on $%
{\cal  {\cal  F}}(W,V),$ and we let $\gamma _{SK}^{0}(W,V)$ denote the
corresponding normed space. We conclude from (\ref{exactint}) that we have
an isometric injection
\begin{equation}
{\cal  I}^{\rm ex}(V,W)\hookrightarrow \gamma _{SK}^{0}(W,V)^{*},
\label{exactint1}
\end{equation}
and in particular, if $W$ is finite-dimensional, we have that
\begin{equation}
{\cal  I}^{\rm ex}(V,W)=\gamma _{SK}^{0}(W,V)^{*}.  \label{exactint2}
\end{equation}

We also have that the exactly integral norm is local.

\proclaim{Proposition}
\label{exactlocal} A linear mapping $\varphi :V\rightarrow W$ is exactly
integral with $\iota ^{\rm ex}(\varphi )\le 1$ if and only if for every finite-dimensional
 subspace $E\subseteq V$ we
have $\iota ^{\rm ex}(\varphi _{|E})\le 1.$
\endproclaim 

\demo{Proof}%
Since ${\cal  I}^{\rm ex}(V,W)$ is an operator ideal, it is clear that $\iota
^{\rm ex}(\varphi _{|E})\le \iota ^{\rm ex}(\varphi )$. On the other hand, for any
finite-dimensional subspace $F$ of $M_{n}$ and any complete contraction $%
\psi :F\rightarrow V$, let $E=\psi (F)$. We have
$$
\nu (\varphi \circ \psi )=\nu (\varphi _{|E}\circ \psi )\le \iota
^{\rm ex}(\varphi _{|E}).
$$
It follows from Theorem \ref{P4.4a} that
\begin{eqnarray*}
\iota ^{\rm ex}(\varphi ) &=&\sup \{\nu (\varphi \circ \psi ):\left\| \psi
:F\rightarrow V~\right\| _{cb}\leq 1{\rm \thinspace , }\enspace F\hbox{ \rm  a matrix
space}\} \\
&\le &\sup \{\iota ^{\rm ex}(\varphi _{|E}):E\subseteq V~\hbox{\rm finite-dimensional}\}.
\end{eqnarray*}
\enddemo %

If $V$ is a matrix space, it is immediate from Theorem \ref{P4.4a} that $%
{\cal  I}(V,W)={\cal  I}^{\rm ex}(V,W).$ This is true more generally. We
recall from Pisier \cite{Pi0} that an operator space $V$ is called %
$1$-{\it exact} if for every finite-dimensional operator space $E$ of $V$ and $%
\varepsilon >0$ there is a matrix space $F$ and a linear isomorphism\break $%
S:E\rightarrow F$ such that $\left\| S\right\| _{cb}\left\| S^{-1}\right\|
_{cb}<1+\varepsilon $. The following result motivated the terminology ``exact
integral.''

\proclaim{Proposition}
\label{Pfin} An operator space $V$ is $1$\/{\rm -}\/exact if and only if for any operator
space $W${\rm ,}
\begin{equation}
{\cal  I}(V,W)={\cal  I}^{\rm ex}(V,W).  \label{1exact}
\end{equation}
\endproclaim 

\demo{Proof}%
If $V$ is 1-exact, then for any finite rank mapping $\psi :W\rightarrow V, $%
$$
\gamma _{SK}(\psi )=\left\| \psi \right\| _{cb},
$$
since if we let $F=\psi (W),$ and we are given $\varepsilon >0,$ we may find a
diagram
$$
\begin{array}{ccccc}
&  & E &  &  \\
&  & \,\,\,\,{\scriptstyle S}\uparrow \downarrow {\scriptstyle S^{-1}} &  &
\\
W & \overset{\psi }{\longrightarrow } & F & \subseteq & V
\end{array}
$$
with $E$ a matrix space and $\max \left\{ \left\| S\right\| _{cb},\left\|
S^{-1}\right\| _{cb}\right\} <1+\varepsilon $. From (\ref{isom}) and
(\ref{exactint1})  it follows that for any linear mapping $\varphi :V\rightarrow W,$ we
have   $\iota (\varphi )=\iota ^{\rm ex}(\varphi ).$

Let us suppose that (\ref{1exact}) holds for all finite-dimensional
subspaces\break $W\subseteq V.$ Then fixing such a subspace, we have a norm
decreasing linear isomorphism (both sides coincide with the vector space $%
W^{*}\otimes V)$%
$$
\theta :\gamma _{SK}^{0}(W,V)\rightarrow {\cal  CB}^{0}(W,V).
$$
But we are given that the adjoint mapping
$$
\theta ^{*}:{\cal  I}(V,W)\rightarrow {\cal  I}^{\rm ex}(V,W)
$$
is isometric, and thus $\theta $ must itself be an isometry. Letting $%
j:W\hookrightarrow V$ be the inclusion mapping, it follows that for any $%
\varepsilon >0,$ we have a matrix space $E$ and a commutative diagram
$$
\begin{array}{ccccc}
&  & E &  &  \\
& {\scriptstyle a}\nearrow &  & \searrow {\scriptstyle b} &  \\
W &  & \overset{j}{\longrightarrow } &  & V\ ,
\end{array}
$$
where $\left\| a\right\| _{cb}\left\| b\right\| _{cb}<1+\varepsilon .$ Thus  
 $W$ is 1-exact  and we conclude $V$ is exact.
\enddemo 

\proclaim{Proposition}
Given operator spaces $V$ and $W$ and a linear mapping $\varphi
:V\rightarrow W${\rm ,} we have that
$$
\iota (\varphi )\le \iota (\varphi ^{*}).
$$
Moreover{\rm ,} $V$ is locally reflexive if and only if we have the isometry
$$
\iota (\varphi )=\iota (\varphi ^{*})
$$
for all operator spaces $W$ and linear mappings $\varphi :V\rightarrow W.$
\endproclaim 

\demo{Proof}%
For any finite-dimensional operator space $E$ we have that
$$
E^{*}\ \check{\otimes} \ V^{*}=(E\hat{\otimes}V)^{*},
$$
whereas the corresponding mapping
$$
E^{*}\hat{\otimes}V^{*}\rightarrow (E\ \check{\otimes} \ V)^{*}
$$
is norm-decreasing. From these we conclude that
\begin{eqnarray*}
\iota (\varphi ) &=&\sup \{\left\| {\rm id}\otimes \varphi :E^{*}\ \check{\otimes} \ %
V\rightarrow E^{*}\hat{\otimes}W\right\| :~E~\hbox{\rm  finite-dimensional}\} \\
&=&\sup \{\left\| ({\rm id}\otimes \varphi )^{*}:(E^{*}\hat{\otimes}%
W)^{*}\rightarrow (E^{*}\ \check{\otimes} \ V)^{*}\right\| :~E~\hbox{\rm  finite-dimensional}\} \\
&\le &\sup \{\left\| {\rm id}\otimes \varphi ^{*}:E\ \check{\otimes} \ W^{*}\rightarrow
E\hat{\otimes}V^{*}\right\| :~E~\hbox{\rm  finite-dimensional}\} \\
&=&\iota (\varphi ^{*}).
\end{eqnarray*}
If $V$ is locally reflexive, then $E\hat{\otimes}V^{*}\rightarrow (E^{*}%
\ \check{\otimes} \ V)^{*}$ is isometric, and the above calculations show that $%
\iota (\varphi )=\iota (\varphi ^{*}).$

If $W$ is a finite-dimensional operator space, we have the isometries
$$
{\cal  CB}(V,W)\cong V^{*}\ \check{\otimes} \ W\cong W\ \check{\otimes} \ V^{*}=%
{\cal  CB}(W^{*},V^{*}),
$$
and from Lemma \ref{P4.2}, we have the isometry
$$
{\cal  I}(W^{*},V^{*})={\cal  N}(W^{*},V^{*}).
$$
Therefore, if for every $\varphi :V\rightarrow W$, we have $\iota (\varphi
)=\iota (\varphi ^{*})$, then from Lemma \ref{P3.1} we have the isometries
$$
{\cal  I}(V,W)={\cal  I}(W^{*},V^{*})={\cal  N}(W^{*},V^{*})={\cal  N%
}(V,W),
$$
and we conclude from Proposition \ref{P4.3a} that $V$ is locally reflexive.%
\enddemo %

\section{The local reflexivity principle for von Neumann preduals}

\proclaim{Theorem}
\label{P5.1} For any $C^{*}$\/{\rm -}\/algebra $A,$ $A^{*}$ is a locally reflexive
operator space{\rm .}
\endproclaim 

\demo{Proof}%
From Proposition 4.4, it suffices to show that we have the isometry
$$
{\cal  I}(A^{*},F)={\cal  N}(A^{*},F)
$$
for all finite-dimensional operator spaces $F.$ Given $\varphi :A^{*}\to F$,
it is trivial from the definition that $\iota (\varphi )\le \nu (\varphi )$.
On the other hand, we have the mappings
$$
{\cal  I}(A^{*},F)\overset{S}{\longrightarrow }{\cal  I}%
^{\rm ex}(F^{*},A^{**})\cong {\cal  I}(F^{*},A^{**})\cong {\cal  N}%
(F^{*},A^{**})\overset{S^{-1}}{\cong }{\cal  N}(A^{*},F),
$$
where $S(\varphi )=\varphi ^{*}$ is contractive (Lemma \ref{P4.4}), and the
second, third and fourth identifications are proved in Lemmas \ref{P4.5},
\ref{P4.2}, and \ref{P3.1}. Therefore, we must have $\nu (\varphi )=\iota
(\varphi )$.
\enddemo %

It was noted in \cite[p.~185]{ER3}, that any subspace of a locally reflexive
operator space is again locally reflexive. In particular, if $R$ is a von
Neumann algebra, then we may identify $R_{*}$ with an operator subspace of $%
R^{*},$ and we obtain the following result.

\proclaim{{C}orollary}
\label{P5.2}For any von Neumann algebra $R,$ the predual $R_{*}$ is a
locally reflexive operator space{\rm .}
\endproclaim

Guided by the classical theory, we may prove stronger versions of
approximation. We begin by recalling ``Helly's lemma'' (see \cite[p.~73]{DF}).

\proclaim{Lemma}
\label{P5.3} Suppose that $E$ is a Banach space and that $L$ is a finite-dimensional 
subspace of $E^{*}.$ Then given
any element $u\in E^{**}$ and $%
\varepsilon >0,$ there exists an element $u_{0}$ of $E$ such that $\left\|
u_{0}\right\| <(1+\varepsilon )\left\| u\right\| $ and
$$
\langle u,h\rangle =\langle u_{0},h\rangle
$$
for all $h\in L.$
\endproclaim

\proclaim{Lemma}
\label{P5.4} Suppose that $V$ is a locally reflexive operator space{\rm ,} and that
$F\subseteq V^{**}$ and $N\subseteq V^{*}$ are finite-dimensional{\rm .} Then for
each $\varepsilon >0$ there exists a mapping $\varphi :F\rightarrow V$ such
that $\left\| \varphi \right\| _{cb}<1+\varepsilon ${\rm ,} and
$$
\langle \varphi (v),f\rangle =\langle v,f\rangle
$$
for all $f\in N.$
\endproclaim

\demo{Proof}%
Local reflexivity implies that we have the isometry (\ref{F1.0}). We can
regard the inclusion mapping $\iota :F\rightarrow V^{**}$ as a contractive
element of $F^{*}\ \check{\otimes} \ V^{**},$ and $L=F\otimes N$ as a finite-dimensional
 subspace of $(F^{*}\
\check{\otimes} \ V)^{*}.$ From Helly's lemma, we can choose an element $\varphi \in F^{*}\ \check{\otimes} \ V$
such that $%
\left\| \varphi \right\| _{cb}<1+\varepsilon ,$ and
$$
\langle \varphi (v),f\rangle =\langle \varphi ,v\otimes f\rangle =\langle
\iota ,v\otimes f\rangle =\langle v,f\rangle
$$
for all $v\in F$ and $f\in N^{*}.$%
\enddemo %

The following result of Pisier (see \cite[Lemma 7.1.4]{Pi1}), is the
analogue of a well-known theorem in Banach space theory.

\proclaim{Lemma}
\label{P5.5} Suppose that $V$ is an operator space{\rm ,} and that $v_{i}\in V,$ $%
f_{i}\in V^{*}$ ($i=1,\ldots ,n)$ are biorthogonal{\rm ,} i{\rm .}e{\rm .,} $%
f_{i}(v_{j})=\delta _{i,j}.$ Then given $\varepsilon >0$ and elements $w_{i}$
such that
\begin{equation}
\sum_{i}\left\| f_{i}\right\| \left\| v_{i}-w_{i}\right\| <\varepsilon ,
\label{E5}
\end{equation}
it follows that there is a complete isomorphism $\varphi :V\rightarrow V$
such that $\varphi (v_{i})=w_{i},$ where $\left\| \varphi \right\| _{cb}\le
1+\varepsilon $ and $\left\| \varphi ^{-1}\right\| _{cb}\le (1-\varepsilon )^{-1}.$
\endproclaim

We say that an operator space $V$ is {\it strongly locally reflexive}, if
given a finite-dimensional subspace $F\subseteq V^{**}$ and a finite-dimensional 
subspace $N\subseteq V^{*},${\it \
}there exists a complete isomorphism $\varphi \ $of $F$ onto a subspace $E$ of $V$ such that

\begin{itemize}
\item[(a)]  $\left\| \varphi \right\| _{cb},\left\| \varphi ^{-1}\right\|
_{cb}<1+\varepsilon ,$

\item[(b)]  $\langle \varphi (v),f\rangle =\langle v,f\rangle $ for all $%
v\in F$ and $f\in N,$ and

\item[(c)]  $\varphi (v)=v$ for all $v\in F\cap V.$
\end{itemize}

We will see in the next section that since the $C^{*}$-algebra $M_{\infty
}=(K_{\infty })^{**}$ contains every finite-dimensional operator space, $%
K_{\infty }$ is not strongly locally reflexive.

\proclaim{Theorem}
\label{P5.6} Suppose that $V$ is a locally reflexive operator space for
which there exists a completely isometric injection
$$
\theta :V^{**}\rightarrow B(H)
$$
which satisfies the ${\rm W}^{*}{\rm AP}$ {\rm (}\/see \S {\rm 2).}  Then $V$ is strongly locally
reflexive{\rm .}
\endproclaim 

\demo{Proof}%
In order to simplify the notation, we will assume that $V^{**}\subseteq
B(H), $ and that we have a net of weak$^{*}$ continuous finite rank complete
contractions $t_{\lambda }:V^{**}\rightarrow B(H)$ for which
$$
\left\| t_{\lambda }(v^{**})-v^{**}\right\| \rightarrow 0
$$
for each $v^{**}\in V^{**}.$

We fix $0<\varepsilon <1$ and finite-dimensional subspaces $F\subseteq V^{**}$ and $N\subseteq V^{*}$, and we let
$0<\delta
\le \varepsilon /3n^{3}$, where $n={\rm dim\thinspace }F$. Then we can choose
a complete contraction $t=t_{\lambda }\ $ such that
\begin{equation}
\left\| t(v)-v\right\| <\delta \left\| v\right\|  \label{E3}
\end{equation}
for all $v\in F.$ Letting $W$ be the range of $t,$ we have $%
t:V^{**}\rightarrow W$ and\break $t^{*}:W^{*}\rightarrow V^{*}.$

From Lemma \ref{P5.4}, there exists a mapping $\varphi :F\rightarrow V$ such
that
\begin{equation}
\left\| \varphi \right\| _{cb}<1+\delta ,  \label{E3a}
\end{equation}
and
\begin{equation}
\langle \varphi (v),f\rangle =\langle v,f\rangle  \label{E4}
\end{equation}
for all $v\in F$ and $f\in t^{*}(W^{*})+N$. We let $C\subseteq {\cal  CB}%
(F,V)$ be the convex set of all mappings $\varphi :F\rightarrow V$
satisfying (\ref{E3a}) and (\ref{E4}). We let $F_{0}=F\cap V$, and $\iota
_{0}:F_{0}\rightarrow V$ be the inclusion mapping. We let $C_{0}\subseteq
{\cal  CB}(F_{0},V)$ denote the convex set of all mappings $\varphi \circ
\iota _{0},$ where $\varphi \in C$. We claim that $\iota _{0}$ is in the
point-norm closure of $C_{0}.$ This is apparent since if we are given an
arbitrary finite-dimensional subspace $G\subseteq V^{*},$ then our previous
argument shows that there is a mapping $\varphi ^{\prime }:F\rightarrow V$
satisfying
$$
\left\| \varphi ^{\prime }\right\| _{cb}<1+\delta ,
$$
and
$$
\langle \varphi ^{\prime }(v),f\rangle =\langle v,f\rangle
$$
for all $v\in F$ and $f\in t^{*}(W^{*})+N\,\,+G$. Since $\iota
_{0}(F_{0})\subseteq V,$ we can suitably choose a net of $\varphi ^{\prime }$
such that $\varphi ^{\prime }\circ \iota _{0}$ converges to $\iota _{0}$ in
the point-weak topology. The usual convexity argument, implies that $\iota
_{0}$ is in the point-norm closure of $C_{0}$, and since $F_{0}$ is finite-dimensional, 
we may choose a mapping
$\varphi \in C\ $ satisfying (\ref{E3a}) and (\ref{E4}), for which
$\Vert \iota_0-\varphi \circ \iota_0\Vert <\delta$ and thus
$$
\left\| \iota _{0}-\varphi \circ \iota _{0}\right\| <n \delta .
$$
For all $v\in F$ and $f\in W^{*}$, we have that
$$
\langle t(\varphi (v)),f\rangle =\langle \varphi (v),t^{*}(f)\rangle
=\langle v,t^{*}(f)\rangle =\langle t(v),f\rangle ;
$$
  thus
\begin{equation}
t(\varphi (v))=t(v)  \label{localeq}
\end{equation}
for all $v\in F.$

We next perturb $\varphi $ in order to satisfy (c). It follows from \cite[Lemma 5.2]{ER4}  that there is a projection
$P$ of $F$ onto $F_{0}=V\cap F$ with $1\leq \left\| P\right\| _{cb}\le n^{2}$. Then
$$
\varphi _{1}=(\iota _{0}-\varphi )P+\varphi :F\rightarrow V
$$
is a completely bounded mapping such that $\varphi _{1}(v_{0})=v_{0}$ for $%
v_{0}\in F_{0},$ and if $v\in F,$%
$$
\langle \varphi P(v),f\rangle =\langle P(v),f\rangle
$$
and thus
$$
\langle \varphi _{1}(v),f\rangle =\langle \varphi (v),f\rangle =\langle
v,f\rangle
$$
for $f\in N,$ i.e., $\varphi _{1}$ satisfies (b) and (c). We also have that
\begin{equation}
\left\| \varphi _{1}\right\| _{cb}\le \left\| \iota _{0}-\varphi \circ \iota
_{0}\right\|_{cb} \left\| P\right\| _{cb}+(1+\delta )\le \delta n^{3}+(1+\delta
)<1+\varepsilon . \qquad \label{E6}
\end{equation}

We let $E$ be the range of $\varphi _{1}.$ We must show that $\varphi _{1}$
is almost a complete isometry of $F$ onto $E$. Let us assume that $%
v_{1},\ldots ,v_{n}$ is an Auerbach basis for $F$ with bi-orthogonal dual
basis $f_{i}$ (i.e., $\left\| v_{i}\right\| =1$, $\left\| f_{i}\right\| =1$
and $f_{i}(v_{j})=\delta _{ij}$). For each $i,$ we have from (\ref{localeq})
that
\begin{eqnarray*}
\left\| v_{i}-t\varphi _{1}(v_{i})\right\| &\le &\left\| v_{i}-t\varphi
(v_{i})\right\| +\left\| t\varphi (v_{i})-t\varphi _{1}(v_{i})\right\| \\
&\leq &\left\| v_{i}-t(v_{i})\right\| +\left\| t\right\| _{cb}\left\|
\varphi \circ \iota _{0}-\iota _{0}\right\| \left\| P\right\| _{cb} \\
&\le &\delta +(1+\delta )\delta \left\| P\right\| _{cb} \\
&\le &\delta +2\delta \left\| P\right\| _{cb}\leq 3\delta \left\| P\right\|
_{cb}<3\delta n^{2}.
\end{eqnarray*}
Thus we have that
$$
\sum_{i}\left\| v_{i}-t\varphi _{1}(v_{i})\right\| \left\| f_{i}\right\| \le
3\delta n^{3}<\varepsilon .
$$
From Lemma \ref{P5.5}, we may find a mapping $s:B(H)\rightarrow B(H)$ for
which
$$
st\varphi _{1}(v_{i})=v_{i}.
$$
and $\left\| s\right\| _{cb}\le (1-\varepsilon )^{-1}.$ It follows that $%
\varphi _{1}^{-1}=st|_{E},$ and since $t$ is completely contractive,
\begin{equation}
\left\| \varphi _{1}^{-1}\right\| _{cb}\le (1-\varepsilon )^{-1}.  \label{E7}
\end{equation}
Then $\varphi _{1}$ will also satisfy (a).%
\enddemo %

\proclaim{Theorem}
\label{P5.7}If $R$ is a von Neumann algebra{\rm ,} then $R_{*}$ is strongly
locally reflexive{\rm .}
\endproclaim 

\demo{Proof}%
\enspace We have from Theorem \ref{P2.1} that any complete isometry\break $R^{*}\rightarrow
B(H)$ has the ${\rm W}^{*}{\rm AP},$ and thus since $R_{*}$ is locally reflexive
(Cor.~\ref{P5.2}), the result follows from Theorem \ref{P5.6}.
\enddemo %

\section{Finite representability and factorizations}

The Banach space notion of {\it finite representability} was introduced by
R.C. James in \cite{James}. It has proved to be quite useful (for example,
see Heinrich \cite{He}), and it has an obvious analogue in operator space
theory.

Given operator spaces $E$ and $F$ and $\varepsilon >0$, we write $E\overset{%
1+\varepsilon }{\cong }F$ if $E$ is $(1+\varepsilon )$-completely isomorphic to $F$; i.e., there is a linear isomorphism
$T:E\rightarrow F$ such that $\left\| T\right\| _{cb}\left\| T^{-1}\right\| _{cb}<1+\varepsilon $. This is equivalent
to saying that the completely bounded Banach-Mazur distance introduced in
\cite{Pi0} satisfies
$$
d_{cb}(E,F)=\inf \left\{ \left\| S\right\| _{cb}\left\| S^{-1}\right\|
_{cb}:S:E\cong F\right\} <1+\varepsilon .
$$
We write $E\overset{1+\varepsilon }{\subseteq }F$ if there is a subspace $%
F_{0}\subseteq F$ with $E\overset{1+\varepsilon }{\cong }F_{0}.$

Let us suppose that we are given a family of operator spaces ${\cal  W}$ .
We say that an operator space $V$ is {\it finitely representable in} $%
{\cal  W}$ if for every finite-dimensional subspace $E$ of $V$ and $%
\varepsilon >0$, $E\overset{1+\varepsilon }{\subseteq }W$ for some $W\in {\cal  W%
}$. If ${\cal  W}=\{W\},$ we simply say that $V$ is finitely representable
in $W$. Two operator spaces $V$ and $W$ are {\it finitely equivalent }if
each is finitely representable in the other. Turning to an important
example,  from the work of Kirchberg and Pisier (\cite{Ki}, \cite{Pi0})  we see that an operator space $V$ is 1-exact
if and only if it is finitely representable in $\left\{ M_{n}\right\} _{n\in \Bbb{N}}$, or equivalently,
it is finitely representable in $K_{\infty }${\it . }The following is an
immediate consequence of Theorem \ref{P5.7}.

\proclaim{{C}orollary}
\label{P6.2}If $R$ is an arbitrary von Neumann algebra{\rm ,} then $R^{*}$ is
finitely equivalent to $R_{*},$ and if $A$ is a $C^{*}$\/{\rm -}\/algebra{\rm ,} $A^{***}$
is finitely equivalent to~$A^{*}.$
\endproclaim

We note that Corollary \ref{P6.2} is false for $C^{*}$-algebras. We have,
for example, that although $K_{\infty }$ is finitely representable in $%
\{M_{n}\}_{n\in \Bbb{N}}$, that is not the case for $M_{\infty }=K_{\infty
}^{**}$.

We may use Corollary \ref{P6.2} to formulate other invariants for the
preduals of von Neumann algebras which are preserved by taking second duals.
To illustrate this, we again have from \cite[Prop.~4.3]{ER4} (or see below)
that for any Hilbert space $H$, $T(H)$ is finitely representable in $%
\{T_{n}\}_{n\in \Bbb{N}}$. Thus an operator space $V$ is finitely
representable in $\{T_{n}\}_{n\in \Bbb{N}}$ if and only if $V$ is finitely
representable in $T_{\infty }$. We have from Corollary \ref{P6.2} that the
predual $R_{*}$ of a von Neumann algebra $R$ is finitely representable in $%
\{T_{n}\}_{n\in \Bbb{N}}$ if and only if that is true for $R^{*}$.

A $C^{*}$-algebra $A$ is said to have the {\it weak expectation property }(%
{\rm WEP}) of Lance \cite{La} if given any faithful representation $%
A\hookrightarrow B(H)$, there is a completely positive contraction $P$ of $%
B(H)$ into the weak closure $\bar{A}$ such that $P(a)=a$ for all $a\in A.$
It is well known that nuclear $C^{*}$-algebras and injective $C^{*}$%
-algebras have the {\rm  WEP}. A $C^{*}$-algebra $A$ has the {\rm WEP} if and only if
given the universal representation $A\subseteq {\cal  B}(H),$ there is a
complete contraction $P:{\cal  B}(H)\rightarrow \bar{A}$ such that $P(a)=a$
for all $a\in A.$ $A$ is said to have the {\rm QWEP} if it is a $C^{*}$%
-algebraic quotient of a {\rm  WEP} algebra. It has been conjectured that
{\it all} $C^{*}$-algebras have the {\rm  QWEP}. Kirchberg \cite{Ki1} has shown
that this problem is equivalent to Connes' question of whether any $II_{1}$
factor on a separable Hilbert space can be realized as a subalgebra of the
ultrapower $N^{\omega },$ where $N$ is the hyperfinite $II_{1}$ factor.

\proclaim{{C}orollary}
\label{P6.3}If $A$ is a {\rm QWEP} $C^{*}$\/{\rm -}\/algebra{\rm ,} then $A^{*}$ is finitely
representable in $\{T_{n}\}_{n\in \Bbb{N}}${\rm .}
\endproclaim

\demo{Proof}%
Let us suppose that $A$ is a unital {\rm WEP} $C^{*}$-algebra, and that\break $%
A\subseteq B(H).$ We let $\rho :B(H)^{*}\rightarrow A^{*}$ denote the
restriction mapping. We have from \cite{EL} Theorem 6.3 (i) that there is a
dilation for $A^{*},$ i.e., a completely positive state preserving mapping $%
\theta :A^{*}\rightarrow B(H)^{*}$ such that $\rho \circ \theta ={\rm id}.$ The
adjoint mapping $\theta ^{*}:B(H)^{**}\rightarrow A^{**}$ is again a
completely positive mapping and it preserves the identity. It follows that $%
\theta ^{*}$ is completely contractive, and thus the same is true for $%
\theta .$ Therefore, $\theta $ is a complete isometry, and we may identify $%
A^{*}$ with a (complemented) subspace of $B(H)^{*}.$ Since $T(H)$ and thus $%
B(H)^{*}$ are finitely representable in $\{T_{n}\}_{n\in \Bbb{N}}$, the same
is true for $A^{*}.$

Given a nonunital {\rm WEP }$C^{*}$-algebra $A$, we have from the
discussion in \cite[pp.~458--459]{Ki1}, that the unital extension $A_{1}$ has
the ${\rm  WEP}.$ It follows that $A_{1}^{*}$ is finitely representable in $%
\{T_{n}\}_{n\in \Bbb{N}}$, and since $A^{*}$ may be identified with a
subspace of $A_{1}^{*},$ the same is true for $A^{*}.$

Finally suppose that $A$ is a {\rm QWEP} $C^{*}$-algebra. If we let $J$ be a
closed ideal in a {\rm WEP} $C^*$-algebra $B$ with $A=B/J,$ we may identify $A^{*}$ with the
annihilator $J^{\bot }\subseteq B^{*}.$ Since $B^{*}$ is finitely
representable in $\{T_{n}\}_{n\in \Bbb{N}}$, the same is true for $A^{*}.$
\enddemo %

\proclaim{{C}orollary}
\label{P6.4} If $R$ is a ${\rm  QWEP}$ von Neumann algebra{\rm ,} then $R_{*}$ is
finitely representable in $\{T_{n}\}_{n\in \Bbb{N}}${\rm .}
\endproclaim

As in the classical case, the theory of finite representability is related
to ultraproducts.

\proclaim{Proposition}
\label{P6.5} Let ${\cal  W}$ be a family of operator spaces{\rm .} If $V$ is
finitely representable in ${\cal  W}${\rm ,} then there exists an index set $I$
and an ultrafilter ${\cal  U}$ on $I$ such that for each $\alpha \in I$
there exists an operator space $W_{\alpha }\in {\cal  W}$ such that $V$ is
completely isometric to a subspace of $\prod\limits_{{\cal  U}}W_{\alpha }$%
{\rm .}
\endproclaim 

\demo{Proof}%
Let $I$ be the collection of all pairs $\alpha =(E,\varepsilon )$ with $E$ a
finite-dimensional subspace of $V$ and $\varepsilon >0$. For convenience, we
write $\alpha =(E_{\alpha },\varepsilon _{\alpha })$. There is a canonical
partial order on $I$ given by $\alpha \preceq \alpha ^{\prime }$ if and
only if $%
E_{\alpha }\subseteq E_{\alpha }^{\prime }~$and $\varepsilon _{\alpha }\geq
\varepsilon _{\alpha }^{\prime }.$ For each $\alpha \in I$, we let $I_{\alpha
}=\{\alpha ^{\prime }\in I:~\alpha \preceq \alpha ^{\prime }\}.$ The
collection ${\cal  I}$ of all such sets $I_{\alpha }$ is a filter on $I$,
and it
is evident that
$$
 \bigcap_{\alpha \in I}I_{\alpha }=\emptyset .
$$
We let ${\cal  U}$ on $I$  be a free ultrafilter containing ${\cal  %
{\cal  I}}$.

For each $\alpha =(E_{\alpha },\varepsilon _{\alpha })\in I$, there exists an
element, say $W_{\alpha }$, in ${\cal  W}$ such that $E_{\alpha }$ is $%
(1+\varepsilon _{\alpha })$-completely bounded isomorphic to a finite-dimensional
 subspace $F_{\alpha }$ of
$W_{\alpha }$. For each such $\alpha ,$ we choose a completely bounded isomorphism $s_{\alpha }:E_{\alpha }\to
F_{\alpha }$ such that $\left\| s_{\alpha }\right\| _{cb}\left\| s_{\alpha
}^{-1}\right\| _{cb}<1+\varepsilon _{\alpha }$. We extend $s_{\alpha }$
(nonlinearly) to $V$ by letting $s_{\alpha }(v)=0$ if $v\notin E_{\alpha }.$
It is easy to verify that the mapping
$$
J:V\to \prod\limits_{{\cal  U}}W_{\alpha };~v\mapsto \widehat{(s_{\alpha
}(v))}
$$
is linear and completely isometric from $V$ into $\prod\limits_{{\cal  U}%
}W_{\alpha }$.%
\enddemo %

In contrast to Banach space theory, the converse of above theorem fails for
operator spaces. Any separable operator space $E$  may be realized as an
operator subspace of $\prod\limits_{{\cal  U}}M_{n}$ for a free
ultrafilter ${\cal  U}$ on $\Bbb{N}$, whereas only exact operator
subspaces of $\prod\limits_{{\cal  U}}M_{n}$ are finitely representable in
$\left\{ M_{n}\right\} _{n\in \Bbb{N}}.$ On the other hand, we have a
necessary and sufficient condition for the operator subspaces of the
ultraproduct of $\{T_{n}\}_{n\in \Bbb{N}}$.

\proclaim{Theorem}
\label{P6.6} An operator space $V$ is finitely representable in $%
\{T_{n}\}_{n\in \Bbb{N}}$ if and only if there exists an index set $I$ and a
free ultrafilter ${\cal  U}$ on $I$ such that $V$ is completely isometric
to a subspace of $\prod\limits_{{\cal  U}}T_{n(\alpha )}${\rm .}
\endproclaim 

\demo{Proof}%
If $V$ is finitely representable in $\{T_{n}\}_{n\in \Bbb{N}}$, then from
Proposition \ref{P6.5} there exists an index set $I$ and a free ultrafilter $%
{\cal  U}$ on $I$ such that $V$ is completely isometric to a subspace of $%
\prod\limits_{{\cal  U}}T_{n(\alpha )}$. Conversely, we have a  sequence of
inclusions
$$
\prod\limits_{{\cal  U}}T_{n(\alpha )}\hookrightarrow \prod\limits_{%
{\cal  U}}M_{n(\alpha )}^{*}\hookrightarrow \ell ^{\infty }(I,M_{n(\alpha
)})^{*}.
$$
Since $\ell ^{\infty }(I,M_{n(\alpha )})$ is an injective von Neumann
algebra, its dual space is finitely representable in $\{T_{n}\}_{n\in \Bbb{N}%
}$, and thus the same for any subspace of $\prod\limits_{{\cal  U}%
}T_{\infty }.$%
\enddemo %

It follows from the above argument  and [21, Prop.~1.3] that the converses of Corollaries 7.2 and 7.3 are also true. 
Thus the {\rm QWEP} conjecture mentioned above is true if and only if the predual of any von Neumann algebra is
finitely representable in $\{T_n\}_{n\in \Bbb N}$.  We conclude with an application of strong local reflexivity to a
factorization theorem. The following result was first demonstrated by the second author \cite{J} in response to a
question posed by G. Pisier. The proof used the Kaplansky density theorem, the fact that the completely
1-summing norm $\pi _{1}$ is in trace duality with the mapping norm $\gamma
_{K}$ defined by factorizations through $K_{\infty },$ and Pisier's
ultraproduct factorization characterization for $\pi _{1}.$ In turn, this
result was used to show that $T(H)$ is locally reflexive, the first instance
of Theorem \ref{P5.7}. Here we proceed in the reverse direction, using the
strong local reflexivity theorem to prove the factorization result.

\proclaim{Theorem}
\label{P5.8} Suppose that $V$ and $W$ are finite-dimensional operator spaces
and that $\varphi :V\rightarrow W$ is a linear mapping{\rm .} If $\varphi $ has a
completely bounded factorization
\begin{eqnarray}
&&\label{cd1}
\\
\noalign{\vskip-16pt}
&& 
\begin{array}{ccccc}
&  & B(H) &  &  \\
& {\scriptstyle r}\nearrow &  & \searrow {\scriptstyle s} &  \\
V &  & \overset{\varphi }{\longrightarrow } &  & W,
\end{array}
\nonumber 
\end{eqnarray}
then for each $\varepsilon >0$, there is a factorization
\begin{eqnarray}
&&\label{cd2}\\
\noalign{\vskip-16pt}
&&
\begin{array}{ccccc}
&  & M_{n} &  &  \\
& {\scriptstyle \tilde{r}}\nearrow &  & \searrow {\scriptstyle \tilde{s}} &
\\
V &  & \overset{\varphi }{\longrightarrow } &  & W
\end{array}
\nonumber
\end{eqnarray}
for some $M_{n}$ such that $\left\| \tilde{r}\right\| _{cb}\left\| \tilde{s}%
\right\| _{cb}<\left\| r\right\| _{cb}\left\| s\right\| _{cb}+\varepsilon ${\rm .}
\endproclaim 

\demo{Proof}%
Assume that we have the commutative diagram of complete contractions (\ref
{cd1}). Taking adjoints, we obtain the completely contractive diagram
$$
\begin{array}{ccccc}
&  & B(H)^{*} &  &  \\
& {\scriptstyle \ s^{*}}\nearrow &  & \searrow {\scriptstyle r^{*}} &  \\
W^{*} &  & \overset{\varphi ^{*}}{\longrightarrow } &  & V^{*}\ .
\end{array}
$$
Then $E=s^{*}(W^{*})$ is a finite-dimensional subspace of $B(H)^{*}$, and $%
F=r(V)$ is a finite-dimensional subspace of $B(H)$. It follows from Theorem
\ref{P5.7} that there exists a completely bounded isomorphism $\psi $ from $%
E $ onto a subspace   $E_{\psi }=\psi (E)$ of $T(H)$ such that $\left\|
\psi \right\| _{cb}\left\| \psi ^{-1}\right\| _{cb}<(1+\varepsilon )^{1/2}$ and
$$
\left\langle \psi (s^{*}(f)),r(v)\right\rangle =\left\langle
s^{*}(f),r(v)\right\rangle
$$
for all $f\in W^{*}$ and $v\in V$.

It is shown in \cite{ER4} (and it follows very easily by truncation from
Lemma \ref{P5.5}) that $T(H)$ is a ${\cal  T}$-space, i.e., for any finite-dimensional
 subspace $G$ of $T(H)$ and
$\varepsilon >0$, there exists an $n\in
\Bbb{N}$ and a subspace $\tilde{G}$ of $T(H)$ containing $G$ such that $%
\tilde{G}\overset{1+\varepsilon }{\cong }T_{n}.$ Applying this to $E_{\psi },$
we have a subspace $\tilde{E}_{\psi }\supseteq E_{\psi }$ and a linear
isomorphism $t:\tilde{E}_{\psi }\rightarrow T_{n}$ for which $\left\|
t\right\| ,\left\| t^{-1}\right\| <(1+\varepsilon )^{1/2}.$ Then $t\circ \psi
\circ s^{*}:W^{*}\rightarrow T_{n}$ and $r^{*}\circ t^{-1}:T_{n}\rightarrow
V^{*}$ are completely bounded maps such that
$$
\left\| t\circ \psi \circ s^{*}\right\| _{cb}\left\| r^{*}\circ
t^{-1}\right\| _{cb}\le \left\| t\right\| _{cb}\left\| \psi \right\|
_{cb}\left\| t^{-1}\right\| _{cb}<1+\varepsilon .
$$
Putting $\tilde{r}=(r^{*}\circ t^{-1})^{*}:V\rightarrow M_{n}$ and $\tilde{s}%
=(t\circ \psi \circ s^{*})^{*}:M_{n}\rightarrow W$, we get
\begin{eqnarray*}
\left\langle f,\tilde{s}\circ \tilde{r}(v)\right\rangle &=&\left\langle
(r^{*}\circ t^{-1})\circ (t\circ \psi \circ s^{*})(f),v\right\rangle \\
&=&\left\langle r^{*}\circ \psi \circ s^{*}(f),v\right\rangle \\
&=&\left\langle \psi \circ s^{*}(f),r(v)\right\rangle \\
&=&\left\langle s^{*}(f),r(v)\right\rangle \\
&=&\left\langle f,s\circ r(v)\right\rangle
\end{eqnarray*}
for all $f\in W^{*}$ and $v\in V$. This shows that $\tilde{s}\circ \tilde{r}%
=s\circ r$ and that $\left\| \tilde{r}\right\| _{cb}\left\| \tilde{s}%
\right\| _{cb}\break <1+\varepsilon $.
\enddemo %

\proclaim{Theorem}
\label{P5.9} Suppose that $V$ and $W$ are finite-dimensional operator spaces
and $A$ is a $C^{*}$\/{\rm -}\/algebra having ${\rm  WEP}${\rm .} If $\varphi :V\rightarrow
W $ has a completely bounded factorization
\begin{eqnarray}
&&\label{cd3}
\\ 
\noalign{\vskip-16pt}
&&\begin{array}{ccccc}
&  & A &  &  \\
& {\scriptstyle r}\nearrow &  & \searrow {\scriptstyle s} &  \\
V &  & \overset{\varphi }{\longrightarrow } &  & W,
\end{array}
\nonumber \end{eqnarray}
then for each $\varepsilon >0${\rm ,} there is a factorization
\begin{eqnarray}
&&\label{cd4}
\\ 
\noalign{\vskip-16pt}
&&
\begin{array}{ccccc}
&  & M_{n} &  &  \\
& {\scriptstyle \tilde{r}}\nearrow &  & \searrow {\scriptstyle \tilde{s}} &
\\
V &  & \overset{\varphi }{\longrightarrow } &  & W
\end{array}
\nonumber \end{eqnarray}
for some $M_{n}$ such that $\left\| \tilde{r}\right\| _{cb}\left\| \tilde{s}%
\right\| _{cb}<\left\| r\right\| _{cb}\left\| s\right\| _{cb}+\varepsilon ${\rm .}
\endproclaim 

\demo{Proof}%
Using the universal representation, we may identify $A^{**}$ with a von
Neumann algebra on a Hilbert space $H$, and we may fix a complete
contraction $P:B(H)\rightarrow A^{**}$ such that $P(a)=a$ for all $a\in A.$
Then we can assume that
$$
\iota _{A}\circ r:V\rightarrow A\hookrightarrow A^{**}\subseteq B(H)
$$
is a completely bounded mapping from $V$ into $B(H)$. Taking the second
adjoint of $s$, we get a completely bounded mapping
$$
s^{**}:A^{**}\rightarrow W^{**}=W
$$
such that $s^{**}\circ \iota _{A}=s$. This gives us a completely bounded
factorization
$$
\begin{array}{ccccc}
&  & B(H) &  &  \\
& {\scriptstyle \iota _{A}\circ r}\nearrow &  & \searrow {\scriptstyle %
s^{**}\circ P} &  \\
V &  & \overset{\varphi }{\longrightarrow } &  & W,
\end{array}
$$
where $\left\| \iota _{A}\circ r\right\| _{cb}=\left\| r\right\| _{cb}$ and $%
\left\| s^{**}\circ P\right\| _{cb}=\left\| s\right\| _{cb}$. Then the
result follows from Theorem \ref{P5.8}.%
\enddemo %

\end{document}